                              \documentclass[11pt]{amsart}  

                    \usepackage{graphicx}    \usepackage{latexsym}

\def\pD{\partial D}  \def\wB{\widetilde{B}} \def\bR{\Bbb{R}}
\def\pE{\partial E}  \def\wS{\widetilde{S}} 
\def\pM{\partial M}  \def\wE{\widetilde{E}} \def\bK{\bar{K}}
\def\pN{\partial N} 
 
\def\pR{\partial R} 

\def\t{\widetilde}  \def\ra{\rightarrow}   

\def\Dgn{\De_{\g \times \{-1\}}} \def\Dmn{\De_{\m \times \{-1\}}} \def\Ddn{\De_{\m' \times \{-1\}}}
\def\Dgp{\De_{\g \times \{ 1\}}} \def\Dmp{\De_{\m \times \{ 1\}}} \def\Ddp{\De_{\m' \times \{ 1\}}}

\def\dgI{{\rm deg}I}
\def\dgP{{\rm deg}P} \def\dg{{\rm deg}} 
\def\dgQ{{\rm deg}Q} 
\def\dgU{{\rm deg}U}
\def\dgW{{\rm deg}W}  

\def\dgWpx{{\rm deg}W_x|_{P_x}}  \def\dgWqx{{\rm deg}W_x|_{Q_x}}
\def\dgWpl{{\rm deg}W_l|_{P_l}}  \def\dgWql{{\rm deg}W_l|_{Q_l}}
\def\dgWpr{{\rm deg}W_r|_{P_r}}  

                     \def\nwpg{\newpage}  
\def\epst{\emptyset}                           \def\nidt{\noindent}
\def\ldq{\l_{\d q}}
\def\prtl{\partial}

\def\Ga{\Gamma} \def\De{\Delta} \def\Si{\Sigma} \def\Th{\Theta}
\def\La{\Lambda} \def\Om{\Omega}

\def\a{\alpha}      \def\b{\beta}   \def\c{\gamma} \def\d{\delta}
\def\e{\varepsilon}    \def\g{\eta}   \def\h{\theta}
      \def\l{\lambda} \def\m{\mu}    
         \def\q{\xi}        
              
\def\z{\zeta}

  \def\C{\mathcal{C}}

\def\pvan{\par\vspace{1mm}\noindent} \def\pvbn{\par\vspace{2mm}\noindent}
\def\pvcn{\par\vspace{3mm}\noindent}
\def\pven{\par\vspace{5mm}\noindent}

\newtheorem{Thm}  {Theorem}     \newtheorem{Prop}{Proposition}[section]  
  \newtheorem{Lem} {Lemma}      [section]    
\newtheorem{Claim}{Claim}       \newtheorem{Cor} {Corollary}  [section]

\def\End{\begin{flushright}$\Box$\end{flushright}}

\def\bgnT{\begin{Thm} \label}   \def\endT{\end{Thm}}   \def\endTR{\end{Thm}\begin{proof}}
\def\bgnP{\begin{Prop} \label}  \def\endP{\end{Prop}}  \def\endPR{\end{Prop}\begin{proof}}
\def\bgnL{\begin{Lem} \label}   \def\endL{\end{Lem}}   \def\endLR{\end{Lem}\begin{proof}}
\def\bgnM{\begin{Claim} \label} \def\endM{\end{Claim}} \def\endMR{\end{Claim}\begin{proof}}
\def\bgnC{\begin{Cor} \label}   \def\endC{\end{Cor}}   \def\endCR{\end{Cor}\begin{proof}}

\def\bgnR{\begin{proof}}                               \def\endR{\end{proof}} 
\def\bgnF{\begin{figure}[htpb!] \centerline}           \def\endF{\end{figure}}
\def\bgnI{\begin{itemize}}                             \def\endI{\end{itemize}}

\def\iclg{\includegraphics} 
\def\captn{\caption}

\title{                The almost alternating diagrams of the trivial knot
}
\author{                                Tatsuya TSUKAMOTO 
}
\address{                   Department of Mathematical Sciences, 
                             School of Science and Engineering,  
                                    Waseda University, 
                                 3-4-1 Okubo Shinjuku-ku, 
                                   Tokyo 169-8555 JAPAN
}    
\thanks{                   email address: tsukamoto@fuji.waseda.jp
} 
\thanks{                  The author acknowledges partial support by 
                         JSPS Research Fellowships for Young Scientists.
}
\date{} 
\begin{document} 
                                 \baselineskip=12pt 
\maketitle
\begin{abstract} 
Bankwitz characterized an alternating diagram representing the trivial knot.
A non-alternating diagram is called almost alternating if
one crossing change makes the diagram alternating.
We characterize an almost alternaing diagram representing the trivial knot.
As a corollary we determine an unknotting number one alternating knot with 
a property that the unknotting operation can be done on its alternating diagram.
\end{abstract}

\section{Introduction}                                                          
\label{S:ITR}                                                                   

              Our concern in this paper is to decide if a given link diagram on $S^2$ 
              represents a trivial link in $S^3$. 
              This basic problem of Knot Theory has been worked in three directions 
              with respect to the properties which we require the diagram to have: 
              closed braid position; positivity; and alternation. 
              We pursue the third direction. For the first direction see \cite{BM-TAMS92} 
              and for the second direction see \cite{Crm-JLMC59} and \cite{St-CM04}.
        \pvcn A link diagram is {\it trivial} if the diagram has no crossings. 
              Obviously a trivial link diagram represents a trivial link.
              A portion of a non-trivial link diagram depicted at the left of 
              Figure \ref{F:24lgs} is called a {\it nugatory crossing}.
              Such a local kink may be eliminated for our purpose.
              Therefore we consider only {\it reduced} link diagrams, i.e.
              link diagrams with no nugatory crossings.
        \pvcn Let $L$ be a link diagram on $S^2$ and let $\hat{L}$ be the link projection 
              obtained from $L$ by changing each crossing to a double point.
              If there is a simple closed curve $C$ on $S^2-\hat{L}$ such that each component of 
              $S^2-C$ contains a component of $\hat{L}$, then we call $L$ {\it disconnected} and
              $C$ a {\it separating curve for} $L$. Otherwise we call $L$ {\it connected}.
        \pvcn A non-trivial link diagram is {\it alternating} if overcrossings and 
              undercrossings alternate while running along the diagram. 
              We know that a reduced alternating link diagram never represents a trivial link.

\bgnT{T:AL}   $(${\rm Crowell} \cite{Crw-ANN59}, {\rm Murasugi} \cite{Mu-JMSJ58}$)$
              A splittable link never admits a connected alternating diagram. \endT

\bgnT{T:AK}   $(${\rm Bankwitz} \cite{Ba-MA30}$)$
              The trivial knot never admits a reduced alternating diagram. \endT
        \pvcn
\bgnF  {\iclg[scale=.9, bb=208 430 576 470]{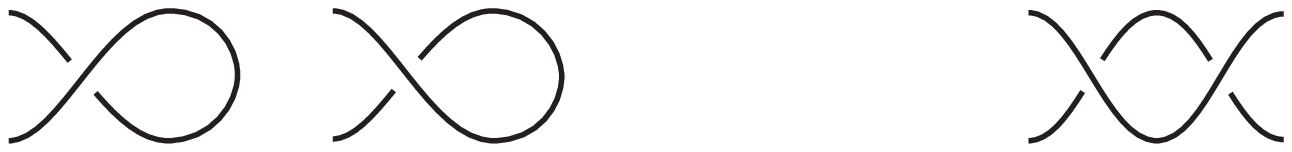}} 
              \captn{Nugatory crossings and a trivial clasp} \label{F:24lgs} \endF
        \pvcn 
              We consider the problem for a link diagram which is alternating except one crossing. 
              Such a link diagram is called almost alternating and 
              first studied by C.Adams et al. \cite{Ad-TPAP92}. 
              A link diagram is {\it almost alternating}  if the diagram is neither trivial 
              nor alternating, but one crossing change makes the diagram alternating. 
              A crossing of an almost alternating link diagram is called a {\it dealternator} 
              if the crossing change at the crossing makes the diagram alternating. 
              In \cite{Ad-TPAP92, Ad-94}, the decision problem  for an almost alternating link
              diagram is asked. M.Hirasawa gave a solution for special almost alternating
              link diagrams in \cite{Hi-TAP00}.
        \pvcn If an almost alternating link diagram has a trivial clasp (the right of 
              Figure \ref{F:24lgs}), then we obtain either a trivial link diagram or 
              an alternating link diagram with fewer crossings from the diagram by the 
              Reidemeister move of type II. Thus we may assume our diagram is
              {\it strongly reduced}, i.e. a reduced diagram with no trivial clasps.
        \pvcn Let $L$ be a non-trivial link diagram on $S^2$ and let $\hat{L}$ be the link 
              projection obtained from $L$ by changing each crossing to a double point.
              If there is a simple closed curve $C$ on $S^2$ intersecting $\hat{L}$ transversely 
              in just two points such that $\hat{L}$ is not a trivial arc in each component of 
              $S^2-C$, then we call $L$ {\it non-prime} and $C$ a {\it decomposing curve for $L$}. 
              Otherwise we call $L$ {\it prime}. Note that a prime link diagram is connected.
        \pvcn We call a portion of a link diagram depicted in Figure \ref{F:6lgs} 
              a {\it flyped tongue}, where the shadowed disks indicate alternating 
              $2$-tangles. Then the author showed the following in \cite{Ts-MPC04}.

\bgnT{T:AAL}  $($\cite{Ts-MPC04}$)$\bgnI      
              \item[$(1)$] A splittable link with $n$-components $(n\geq 3)$ 
                           never admits a connected almost alternating diagram.  
              \item[$(2)$] A prime, strongly reduced almost alternating diagram of 
                           a splittable link with $2$-components has a flyped tongue. 
              \endI \endT    
         \nidt
              The following is the main theorem of this paper.

\bgnT{T:Main}  A strongly reduced almost alternating diagram of the trivial knot has a flyped tongue.
               \endT    
    
\bgnF  {\iclg[scale=.9, bb=77 396 500 466]{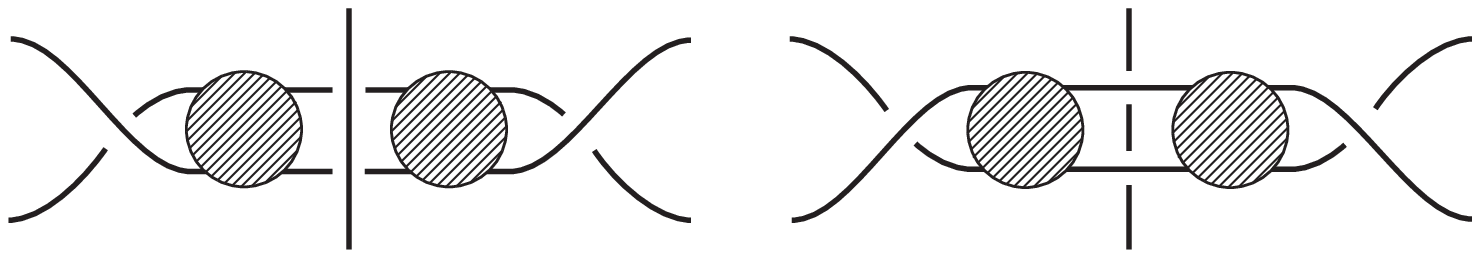}} \captn{Flyped tongues} \label{F:6lgs} \endF

\pven       
              \subsection{The almost alternating diagrams of the trivial knot}
\label{SS:11}
              Theorem \ref{T:Main} yields a simple finite algorithm to see if a given reduced 
              almost alternating knot diagram represents the trivial knot without increasing 
              the number of crossings of diagrams in the process. In fact Adams et.al. in 
              \cite{Ad-93} conjectured that we have a calculs to reduce a given almost alternating 
              diagram of the trivial knot consisting of three kinds of local moves on link diagrams: 
               a {\it flype move} defined by Figure \ref{F:flype}; 
              an {\it untongue move} defined by Figure \ref{F:utgmv}; and
              an {\it untwirl move} defined by Figure \ref{F:utwmv},
              where we allow the move obtained by changing all the crossings in
              or taking the mirror image of each figure.
              Note that each move does not change the link type which a diagram represents.
              The last two moves are introduced in \cite{Ad-93}. We show their conjecture is true.
              A similar algorithm for a reduced almost alternating link diagram with more than
              one component is obtained in \cite{Ts-MPC04}.
        \pven
              Let $K$ be a reduced almost alternating knot diagram. 
              If  $K$ is not strongly reduced, then apply the Reidemeister move of type II 
              to  $K$ to have another diagram $K'$, which is trivial or alternating. 
              In the first case, we can see that $K$ is not reduced, 
              which contradicts the assumption. Consider the second case. 
              Since $K'$ has at most two nugatory crossings, $K'$ represents the 
              trivial knot if and only if $K'$ is a coiled diagram for a non-zero 
              integer $m$ from Theorem \ref{T:AK} (Figure \ref{F:coil}). 
              Next consider the case when $K$ is strongly reduced. If $K$ has no flyped tongues, 
              then $K$ represents a non-trivial knot from Theorem \ref{T:Main}. 
              Otherwise, we obtain another almost alternating diagram $K''$ which has 
              fewer crossings than $K$ by an untongue move or an untwirl move after sufficient 
              flype moves. It is easy to see that $K''$ is reduced, since $K$ is strongly reduced. 
              Then we go back to the beginning and continue the process.
        \pven
              Going over the above process assuming that $K$ represents the trivial knot,
              we obain the following. Here note that  if $K$ is not strongly reduced,
              then $K$ is a diagram in Figure \ref{F:Cm}, which we denote by $\C_m$.

\bgnT{T:Red}  Let $K$ be a reduced almost alternating diagram of the trivial knot. Then there 
              are a non-zero integer $m$ and a sequence of reduced almost alternating diagrams 
              $$K=K_1 \ra \cdots \ra K_p = C_m$$ such that $K_{i+1}$ is obtained from $K_i$ 
              by a flype move, an untongue move or an untwirl move. 
              \endT \End
        \nidt
              Therefore we can obtain all the almost alternating diagrams of the trivial knot.
              Here we defince a {\it tongue move} and a {\it twirl move} as the converse of 
              an untongue move and an untwirl move, respectively. 

\bgnC{C:Gen}  A reduced almost alternating diagram of the trivial knot is obtained from $C_m$ 
              for a non-zero interger $m$ by tongue moves, twirl moves and flype moves. 
              \endC \End
\pvan
\bgnF  {\iclg[scale=.9, bb=108 482 468 525]{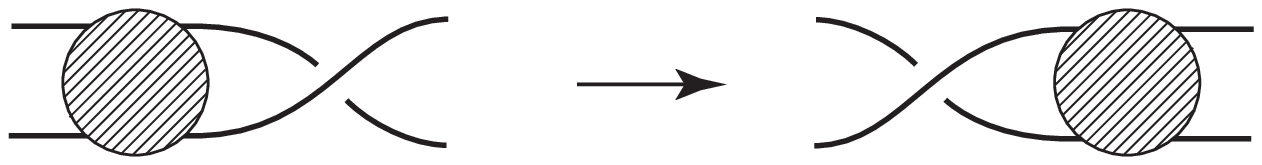}}  \captn{A flype move}    \label{F:flype}\endF
\pvan
\bgnF  {\iclg[scale=.9, bb=110 401 470 461]{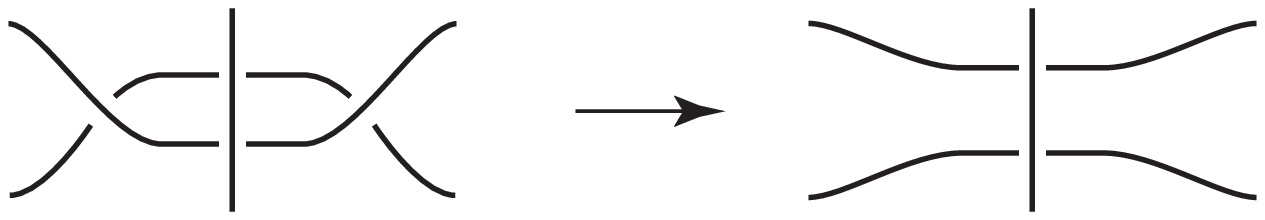}} \captn{An untongue move}\label{F:utgmv}\endF
\pvan
\bgnF  {\iclg[scale=.9, bb=139 427 469 509]{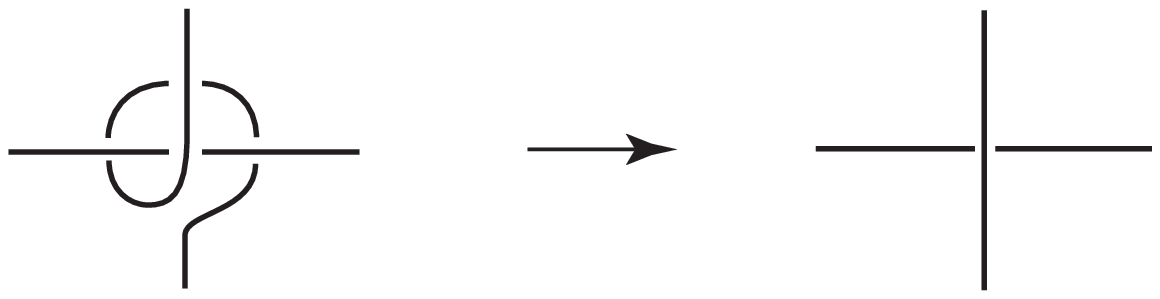}} \captn{An untwirl move}  \label{F:utwmv}\endF
\pvan
\bgnF  {\iclg[scale=.9, bb= 68 419 423 493]{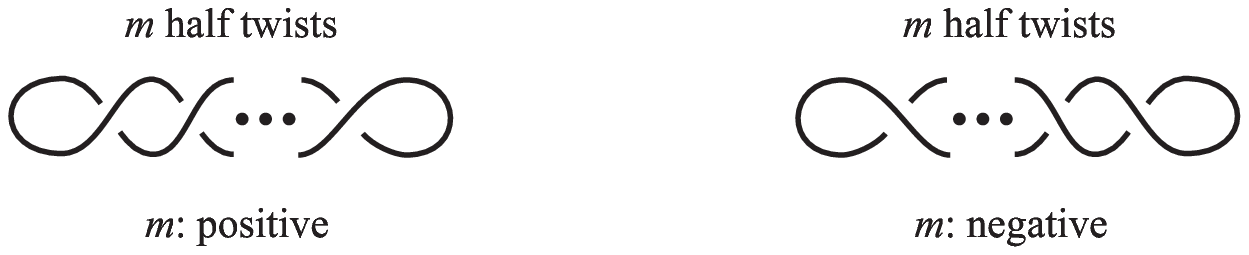}}\captn{Coiled diagrams} \label{F:coil}\endF
\pvan
\bgnF  {\iclg[scale=.9, bb= 66 356 439 493]{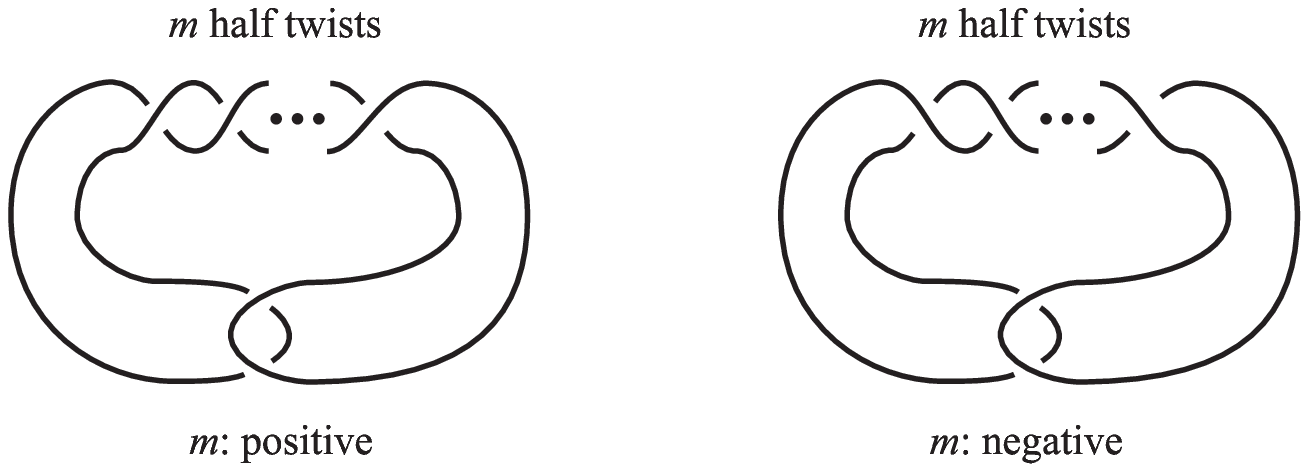}}\captn{$C_m$} \label{F:Cm}\endF
              \subsection{Alternating knots with unknotting number one}
\label{SS:12}
              In \cite{Ko-PAMS91} P.Kohn made a conjecture, which says that a link with
              unknotting number one has a minimal diagram which has a crossing such that the 
              crossing change at the crossing makes the link trivial. This conjecture was 
              shown to be true for large algebraic alternating knots by C.Gordon and J.Luecke 
              in \cite{GL-06}. We remark here that we can obtain all the alternating knots with 
              unknotting number one satisfying the conjecture from Corollary \ref{C:Gen}, since
              we obtain an alternating knot with unknotting number one from a reduced almost 
              alternating diagram of the trivial knot by the crossing change at the dealternator. 
%
%
              \subsection{Organization of the paper}
\label{SS:13}
              Theorem \ref{T:AL} and Theorem \ref{T:AK} were proved algebraically 
              using the Alexander polynomial of a link in \cite{Crw-ANN59} and \cite{Mu-JMSJ58}
              and using the determinant of a knot in \cite{Ba-MA30}, respectively. After those, 
              geometric proofs were given in \cite{WM-TPL84} and in \cite{MT-MPC91} using the  
              ``crossing-ball" technique invented by W.Menasco. Namely he embed a link in a
              ``branched" sphere $S$ to realize a diagram as a geometrical object. We succeed his
              technique to prove Theorem \ref{T:Main} and review it in Section \ref{S:PRE}.
              In Section \ref{S:SPP}, we introduce key concepts which play important roles 
              in this paper: special position for a spanning surface of a link; short arcs; 
              and short bridges. We show that if a given spanning surface is in special position, 
              then the boundary of a neighborhood of it is in standard position
              (Proposition \ref{P:SP/ST}).
              In \cite{Ts-MPC04} to prove Theorem \ref{T:AAL} the author in fact showed that 
              if a prime, strongly reduced almost alternating link diagram on $S$ admits 
              a sphere in its complement which is in standard position, then the diagram 
              admits a flyped tongue (Theorem \ref{T:MPC}). Therefore we are done if our 
              almost alternating diagram on $S$ of the trivial knot admits a spanning disk 
              in special position.
              In Section \ref{S:AAD}, for a spanning surface $E$ of a link which is given 
              as a connected, reduced almost alternating diagram on $S$, we show that $E$ 
              is in special position if and only if $E$ has no short arcs.
              In Section \ref{S:CUT}, we show that if a spanning disk $E$ of the trivial 
              knot which is given as a strongly reduced almost alternating diagram on $S$
              has a short arc, then we can cut $E$ along the short arc or short bridges to 
              have a connected, strongly reduced almost alternating diagram on $S$ of 
              the trivial $2$-component link with spanning disks in special position.
              Then we study the intersection diagram of the spanning disks and $S$ 
              to show that the given diagram has a flyped tongue in Section \ref{S:GLU}.

\pven                                                                           
\section{Preliminary}                                                           
\label{S:PRE}                                                                   

              In this section we bliefly review concepts introduced by Menasco with some 
              additional or modified notations. For details, see \cite{WM-TPL84}, 
              \cite{MT-MPC91} etc.
         \pvan
              Let $\wS$ be a $2$-sphere in $S^3=\bR^3 \cup \infty$. Denote by $\wB^-$ the $3$-ball 
              which $\wS$ bounds in $\bR^3$ and by $\wB^+$ $S^3 - {\rm int}\wB^-$. Take $m$ halls 
              out from $\wS$ and denote the result by $\wS_m$. To each hall, put a $2$-sphere 
              $\h_i$ with an equater $\e_i$ specified so that the equater is on the hall. 
              We call each $\h_i$ a {\it bubble} and the $3$-ball which a bubble bounds in $\bR^3$ 
              a {\it crossing-ball}, denoted by $\Th_i$. We call the disk 
              $\h_i \cap \wB^\pm$ an {\it upper/lower} hemisphere and denote it by $\h_i^\pm$. 
              A {\it bubbled sphere} $S_m$ is a union of $\wS_m$ and the $m$ bubbles.
              We denote the $2$-sphere $\wS_m \cup (\cup \h_i^\pm)$ by $S_m^\pm$ and the 
              $3$-ball which $S_m^\pm$ bounds in $\wB^\pm$ by $B_m^\pm$.
              A link $L$ in $S^3$ is called a {\it link diagram on} $S_m$ if $L$ is on $S_m$,
              meets a bubble $\h_i$ in a pair of two arcs $a_i^+b_i^+$ on $\h_i^+$ and $a_i^-b_i^-$ 
              on $\h_i^-$, and meets the equater transversely so that $a_i^+$, $a_i^-$, $b_i^+$ 
              and $b_i^-$ are positioned on $\e_i$ in this order. 
              Note that $L \cap S_m^+$ on $S_m^+$ is a link diagram in a usual sense. 
              We call a diagram on $S_m$ simply a diagram unless any confusion is expected.
              We say that a link diagram $L$ on $S_m$ has a specific property, e.g.
              alternation, if $L \cap S_m^+$ on $S_m^+$ has the property. 
              Then we also say that $L$ is in alternating {\it position}. 
              We assume that $m$ is sufficiently large and omit $m$ from now on.
        \pvcn
              Let $L$ be an $n$-component link diagram $L_1 \cup \cdots \cup L_n$ on $S$. We call 
              the intersection $L \cap \h_i$ a {\it crossing} if it is 
              not empty. A {\it segment} $\l_j$ is a component of $L \cap (S^+ \cap S^-)$, and a 
              {\it positive/negative long segment} $\La_k^\pm$ is a component of $L \cap S^\pm$. 
              We say that $\La_k^\pm$ {\it runs through} a bubble $\h_i$ if $\La_k^\pm$ contains 
              the arc $a_i^\pm b_i^\pm$ of $L \cap \h_i$, and that $\La_k^\pm$ is $p$-/$n$-{\it 
              adjacent to} $\h_i$ if an end of $\La_k^\pm$ is on $\h_i$. The {\it length} of a 
              long segment $\La_k$ is the number of segments which $\La_k$ has. 
              A segment $\l_j$ is $p$-/$n$-{\it adjacent} to $\h_i$ if the positive/negative long 
              segment containing $\l_j$ is $p$-/$n$-adjacent to $\h_i$. If $\l_j$ is $p$-adjacent 
              and $n$-adjacent to bubbles, then it is called {\it alternating}. Otherwise $\l_j$ 
              is called {\it non-alternating}. A bubble $\h_i$ is $p$-/$n$-{\it adjacent} to 
              another bubble $\h_l$ if there is a segment which has its ends on $\h_i$ and $\h_l$
              and is $p$-/$n$-adjacent to $\h_l$. A crossing $x$ is $p$-/$n$-{\it adjacent} to 
              another crossing $y$ if the bubble at $x$ is $p$-/$n$-adjacent to the bubble at $y$.
              A {\it region} $R_k$ is the closure of a component of $(S^+ \cap S^-) - L$ and  
              its degree, denoted by ${\rm deg}R_k$, is the number of segments on its boundary. 
              Let $N_j$ be a sufficiently small tubular neighborhood of $L_j$ such that
              $\pN_j \cap \Th_i$ is a pair of a saddle-shaped disk in $B^+$ and
              a saddle-shaped disk in $B^-$.


              \subsection{Standard position for a closed surface in a link complement}
\label{SS:21}
              Let $L$ be a link diagram on $S$ and let $F$ be a closed surface in $S^3 - L$. 
              Then we may isotop $F$ so that $F$ satisfies the following conditions,
              and then we say that $F$ is in {\it basic position}. 
              \bgnI
              \item[(Fb1)] $F$ intersects $S^{\pm}$ transversely in a pairwise 
                           disjoint collection of simple closed curves; 
              \item[(Fb2)] $F$ does not intersect $N_j$ for any $j$; and
              \item[(Fb3)] $F$ intersects each crossing-ball $\Th_i$ in a collection of 
                           saddle-shaped disks in $\Th_i-\cup N_j$ $($Figure $\ref{F:XBS}$$)$. 
              \endI
\pvcn     
{\bf Definition}. Let $L$ be a link diagram on $S$ and let $F$ be a closed surface in 
                  $S^3 - L$ which is in basic position with $F \cap S^\pm \neq \epst$. 
                  Let $C$ be a curve of $F \cap S^\pm$. We say that $C$ is {\it standard} 
                  if $C$ satisfies the following conditions and that $F$ is in {\it standard 
                  position} if any curve of $F \cap S^\pm$ is standard.   \bgnI
               \item[(Ft1)] $C$ bounds a disk in $F \cap B^\pm$;
               \item[(Ft2)] $C$ meets at least one bubble; and
               \item[(Ft3)] $C$ meets a bubble in an arc.
               \endI

\bgnF   {\iclg[scale=.9, bb=194 351 315 442]{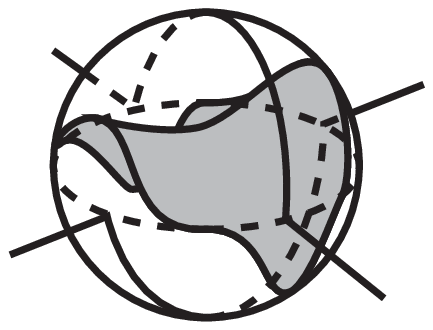}}
         \captn{A saddle-intersection in a crossing-ball}\label{F:XBS}\endF
 
               \subsection{Standard position for a spanning surface of a link}
\label{SS:22}
               Let $E$ be a spanning surface of a link diagram $L$ on $S$. 
               We may isotop $E$ so that $E$ satisfies the following conditions,
               and then we say that $E$ is in {\it basic position}. 
               \bgnI 
               \item[(Eb1)] $E$ intersects $S^\pm$ transversely in a pairwise 
                            disjoint collection of simple closed curves; 
               \item[(Eb2)] $E$ intersects $N_j$ in an annulus $M_j$ so that 
                            $M_j \cap \Th_i = L_j \cap \h_i$ and $\pM_j \cap \pN_j$ 
                            proceeds along $\pN_j$ monotonely with respect to the 
                            longitudinal coordinate of $\pN_j$; and 
               \item[(Eb3)] $E-K$ intersects each crossing-ball $\Th_i$ 
                            in a collection of saddle-shaped disks in $\Th_i-\cup N_j$. \endI
         \pvcn
               We call a component of $(E-K) \cap \h_i^\pm$ a {\it positive/negative}  
               {\it saddle-arc}, and call a component of $\pE \cap S^\pm$ a 
               {\it positive/negative} {\it boundary-arc}. Each end of a boundary-arc is 
               called a {\it junction}. Note that each alternating and non-alternating 
               segment has odd and even number of junctions, respectively.
               The closure of a component of the intersection of $E$ and the interior of 
               a region is an {\it inside arc} if its both ends are on bubbles; an {\it 
               outside arc} if one end is on a bubble and the other is a junction; and 
               an {\it isolated arc} if both ends are junctions. 
               Let $C$ be a curve of $L \cap S^\pm$. We say that $C$ 
               {\it runs through the center} (resp. {\it through a side}) {\it of} $\h_i^\pm$ if 
               $C$ meets $\h_i$ in $L \cap \h_i^\pm$ (resp. in a saddle-arc on $\h_i^\pm$). 
               We say that $C$ {\it runs through} (resp. {\it touches}) a segment $\l_i$ if $C$ 
               meets $\l_i$ in a boundary arc in the interior of $\l_i$ whose end points belong 
               to outside or isolated arcs in different regions (resp. in a same region).
\pvcn
{\bf Definition}. Let $E$ be a spanning surface in basic position of a link diagram $L$ 
               on $S$ and $C$ be a curve of $E \cap S^\pm$. We say that $C$ is {\it standard} 
               if $C$ satisfies the following conditions and that $E$ is in {\it standard 
               position} if any curve of $E \cap S^\pm$ is standard.  \bgnI
               \item[(Et1)] $C$ bounds a disk in $E \cap B^\pm$;
               \item[(Et2)] $C$ meets a bubble or a segment;
               \item[(Et3)] $C$ meets a bubble in an arc.
               \item[(Et4)] $C$ never runs through a side of an upper/lower hemisphere with
                            meeting a segment which is adjacent to the bubble; and
               \item[(Et5)] $C$ never touches a segment. \endI
         \pven
               

               \subsection{Band moves along bridges of a spanning surface of a link}
\label{SS:23}
               Let $E$ be a spanning surface in basic position of a link diagram $L$ on $S$. 
               Assume that there is a disk $\De_\g$ in $B^\pm$ such that: 
               \bgnI
               \item[(Br1)] $\De_\g \cap E     =\g$ is an arc in $\prtl \De_\g$;
               \item[(Br2)] $\De_\g \cap S^\pm =\z$ is an arc in $\prtl \De_\g$;
               \item[(Br3)] $\g \cup \z =\prtl \De_\g$ and $\g \cap \z =\prtl \g=\prtl \z$; and
               \item[(Br4)] $\z$ is in a region $R$.
               \endI
         \pvcn
               We call $\g$ a {\it bridge of} $E$. 
               We say that $\g$ is {\it trivial} if there is a disk $\De'$ in $E \cap B^\pm$
               such that $\prtl \De' = \g \cup \z'$ with $\z'$ in $R$.
               A {\it band move along a bridge} $\g$ is an isotopy performed by sliding
                $\g$ across $\De_\g$ and past $\z$ (see Figure \ref{F:bdmv}). 

\bgnF   {\iclg[scale=.9, bb=99 404 378 594]{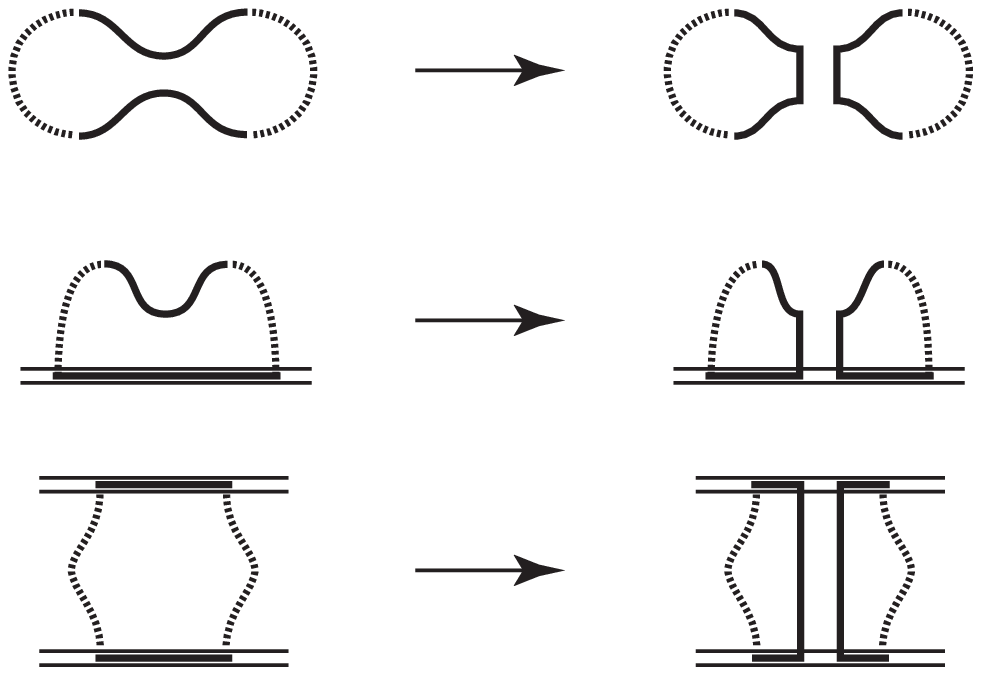}}
         \captn{Three kinds of bridges and the band moves along the bridges} \label{F:bdmv}\endF


               \subsection{The complexity of a spanning surface of a link}
\label{SS:24}
               Let $E$ be a spanning surface in basic position of a link diagram $L$ on $S$. 
               In general $E$ is not in standard position. However $E$ can be isotoped 
               into standard position if $E$ is incompressible.
               Here we define the {\it complexity} of $E$ as the ordered pair $(t, u)$, 
               where $t$ is the number of saddle-intersections of $E \cap \cup \Th_i$ 
               and $u$ is the total number of components of $E \cap S^\pm$. 

\bgnP{P:MC/ST} $($\cite{MT-JRAM92} {\bf Proposition 2.2, 2.3}$)$
               If $E$ is incompressible and has a minimal complexity, 
               then $E$ is in standard position.
               \endP

\section{A spanning surface of a link diagram on $S$}                           
\label{S:SPP}                                                                   

               Let $E$ be an incompressible spanning surface in standard position 
               of a link diagram $L$ on $S$. 
               \subsection{Special position for $E$}
\label{SS:sp}
               Let $C$ be a curve of $E \cap S^\pm$. We say that $C$ is {\it special} 
               if $C$ satisfies the following conditions, and that $E$ is in {\it special 
               position} if any curve of $E \cap S^\pm$ is special.
               \bgnI
               \item[(Ep1)] $C$ never runs through the center of an upper/lower hemisphere with
                            meeting a segment which is adjacent to the bubble;
               \item[(Ep2)] $C$ shares at most $1$ junction with an alternating segment; and
               \item[(Ep3)] $C$ shares no junctions with a non-alternating segment. 
               \endI 
\pven
               Then we have the following.
\bgnP{P:SP/ST} If $E$ is in special position, then the boundary of a neighborhood of $E$ 
               is in standard position.      
               \endPR 
               Let $M'$  be a product neighborhood $E \times [1,-1]$ of $E$ which is 
               sufficiently small compare to the tubular neighborhood $\cup N_j$ of $L$. 
               Take a neighborhood $M$ of $E$ as the union of $M'$ and a neighborhood of $\cup N_j$. 
               Clearly we have that the boundary $\pM$ of $M$ is in basic position and that 
               $\pM \cap S^\pm \neq \epst$. We show only that the positive curves are 
               standard, since we can similarly show that the negative curves are standard.
\pven          It is easy to see that $M \cap S^+$ is a neighborhood of the union of positive 
               curves and positive long segments. Note that a positive long segment of length 
               $p$ meets exactly one positive curve if $p \geq 2$ and no positive curves if 
               $p=1$ from conditions (Ep2) and (Ep3). Therefore we have that $$\pM \cap S^+ 
               = \{C_1', C_1'', \cdots, C_m', C_m'', C_{m+1}, \cdots, C_{m+q}\},$$ where 
               $ E  \cap S^+ = \{C_1, \cdots, C_m\}$ and $C_i'\cup C_i''$ is the boundary 
               of a neighborhood $M_i$ of the union of $C_i$ and the positive long segments 
               which $C_i$ meets, and $C_{m+k}$ is the boundary of a neighborhood of a positive  
               long segment $\La_k$ of length $1$.  
\pven          Then it is clear that $C_{m+k}$ is standard and that $C_i'$ and $C_i''$ 
               satisfy condition (Ft1) from the construction. 
               Note that $C_i$ meets a bubble from conditions (Et2), (Ep2) and (Ep3).
               Therefore $C_i'$ and $C_i''$ satisfy condition (Ft2). 
               Assume that $C_i'$ or $C_i''$, say $C_i'$ does not satisfy condition (Ft3). 
               Note that the pair of the curves of $\pM \cap S^+$ which is closest to the 
               center of an upper hemisphere $\h_k^+$ is the boundary of $M_l$ such that 
               $C_l$ runs through the center of $\h_k^+$. 
               Thus $C_i'$ runs through one side of a bubble twice. 
               This implies that $C_i$ does not satisfy condition (Et3), (Et4) or (Ep1), 
               which is a contradiction (see Figure \ref{F:Ft3}). 
               \endR
\bgnF   {\iclg[scale=.9, bb=90 392 496 479]{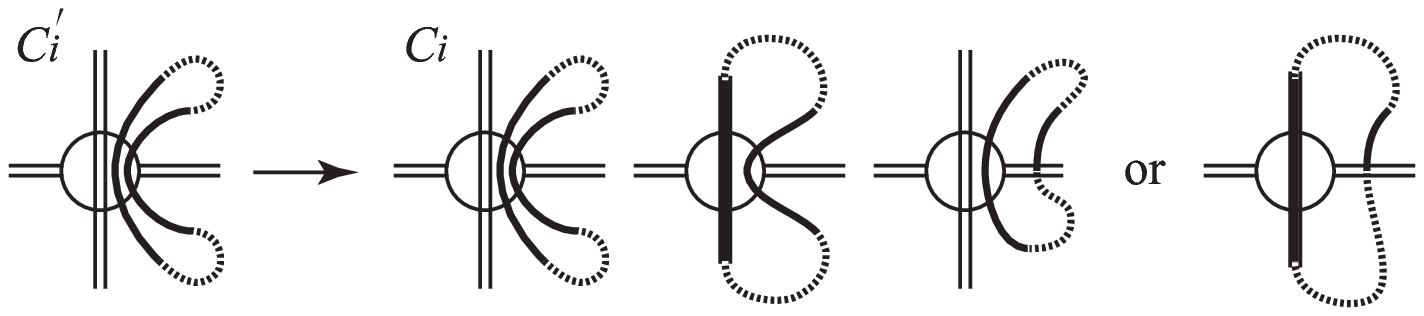}}     \captn{} \label{F:Ft3}\endF

               \subsection{short arcs of $E$}
\label{SS:sa}  A {\it short arc of} $E$ is an isolated arc $\q$ whose ends are on distinct
               segments which are adjacent to a common crossing. Depending upon how the 
               positive curve containing $\q$ meets the segments, we have four types of 
               short arcs as in Figure \ref{F:sarc}, where taking the mirror images do not 
               change their types. The {\it cut surgery along a short arc} $\q$ is the 
               operation of replacing $E$ with $E^\q=E- \q \times (-1,1)$, where we isotop
               $E^\q$ so that $L^\q=\prtl E^\q$ be a link diagram on $S$.
\pvan
\bgnF   {\iclg[scale=.9, bb=94 287 505 480]{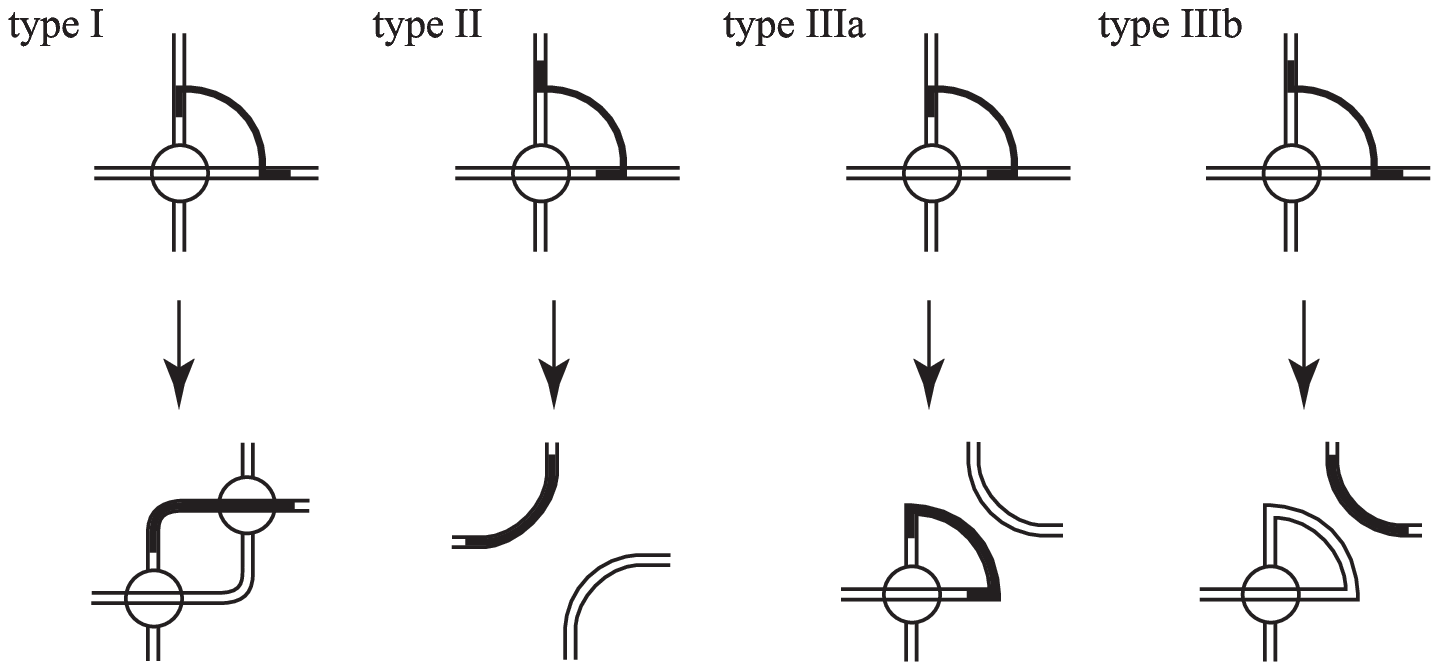}} \captn{} \label{F:sarc}\endF 
%

               \subsection{short bridges of $E$}
\label{SS:sb}  If $\g$ is a non-trivial bridge with its ends on distinct segments which are 
               adjacent to a common crossing $x$, then we call $\g$ a {\it short bridge of} $E$.
               The {\it cut surgery along a short bridge} $\g$ is the operation of replacing $E$ 
               with $E^\g=(E-\g \times (-1,1)) \cup \Dgn \cup \Dgp$ (see Figure \ref{F:sbrdg}). 
\pvan
\bgnF  {\iclg[scale=.9, bb=54 429 456 533]{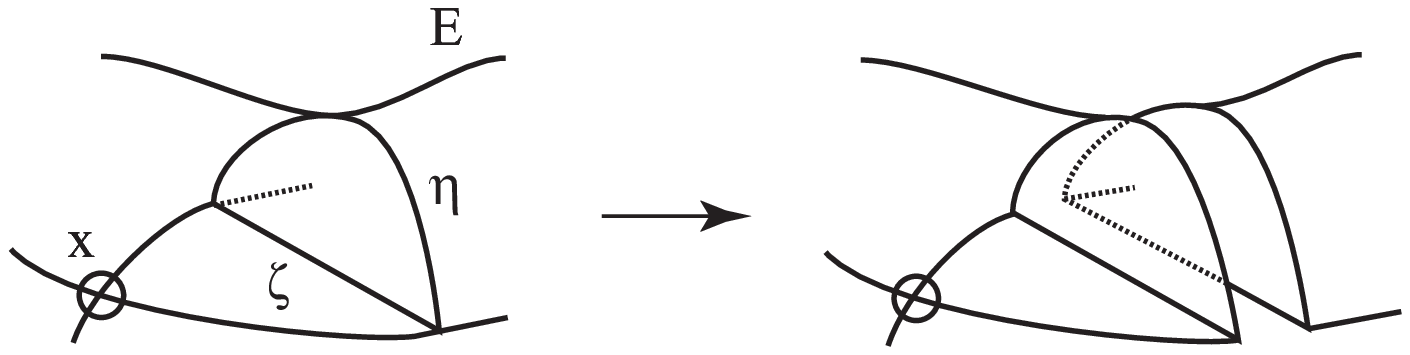}}\captn{}\label{F:sbrdg}\endF

\bgnL{L:3sab} If a crossing $x$ admits a short arc or a short bridge, 
              then $\h_x$ has no saddle-intersections.
              \endLR 
              Assume otherwise. Then, there exists a curve which does not satisfies
              condition (Et3) or (Et4).
              \endR 
%

               \subsection{short cuts of $E$}
\label{SS:sc}  A {\it short cut} $\m$ of $E$ is a short arc of type III or a short bridge.
               The {\it cut surgery along a short cut} $\m$ is the operation of replacing $E$ 
               with $E^\m=(E-\m \times (-1,1)) \cup \Dmn \cup \Dmp$ and
               we let $L^\m=\prtl E^\m$. Note that this is equivalent
               to the cut surgery along a short arc (resp. a short bridge) if $\m$ is a short arc 
               (resp. a short bridge).
        \pvbn
               If a curve does not satisfies condition (Ep1) at a bubble $\h_x$, then we say 
               that the curve {\it has a neck} (at crossing $x$). Then we have the following.

\bgnL{L:3neck} If a curve has a neck at a crossing $x$, then the curve admits a short cut 
               on a region with $x$, and a short arc of type {\rm II} or a non-trivial bridge 
               on another region with $x$.
               \End \endL 
         \pvan 
               Here we define two types of curves each of which consists of two short arcs
               and two boundary arcs: a curve of type $\Ga_1^\pm$ is a curve with one neck 
               which admits a short arc of type II and a short arc of type III; and a curve 
               of type $\Ga_2^\pm$ is a curve with two necks around a non-trivial clasp 
               each of which admits a short arc of type II (see Figure \ref{F:fltc} for 
               curves of type $\Ga_1^+$ and $\Ga_2^+$). 
\pvan
\bgnF   {\iclg[scale=.9, bb=117 431 486 527]{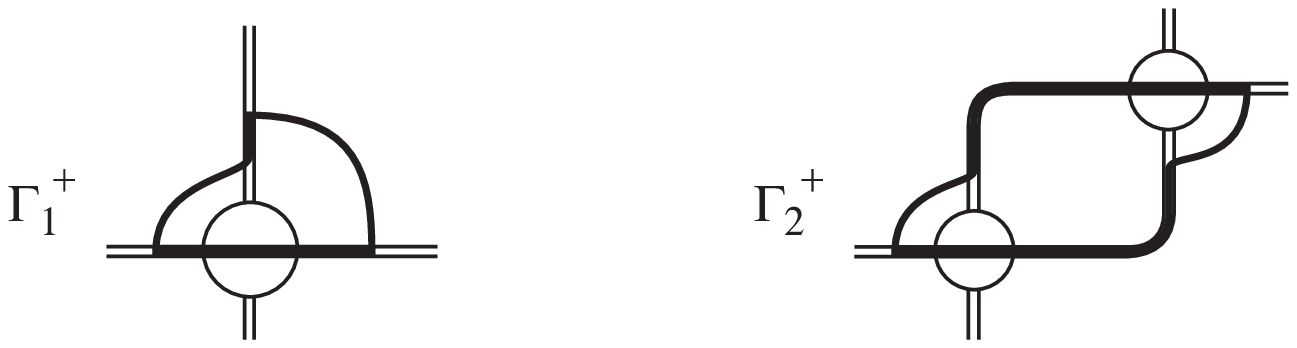}}\captn{} \label{F:fltc}\endF
        \pvan
               Assume that $E$ admits two short cuts $\m$ and $\m'$. We say that $\m$ and $\m'$ 
               are {\it equivalent} if they have the ends on same segments. If $\m$ and $\m'$
               are disjoint and are not equivalent, then let $\Om$ be the subdisk of $E$ bounded 
               by $\m$ and $\m'$. Assuming that $\m \times \{1\}$ and $\m' \times \{-1\}$ belong 
               to $\Om$, define $(\Om^\m)^{\m'}$ as 
               $(\Om^\m)^{\m'}=(\Om - \m \times (0,1) - \m' \times (-1,0)) \cup \Dmp \cup \Ddn.$ 
               We say that $\m$ and $\m'$ are {\it parallel} if each intersection 
               curve of $(\Om^\m)^{\m'} \cap S^\pm$ has type $\Ga_1$ or $\Ga_2$.

\pven                                                                           
\section{A spanning surface of an almost alternating link diagram on $S$}       
\label{S:AAD}                                                                   

              In this section we study an incompressible spanning surface $E$ of a connected, 
              reduced almost alternating link diagram $L$ on $S$. Here we assume that $L$ is 
              not the diagram of Figure \ref{F:2dalt}. 
              Thus $L$ has only one dealternator from 
              the following proposition, and we denote the dealternator by $\d$.
              We denote the bubble at a crossing $x$ by $\h_x$.
              We call a curve of $E \cap S^\pm$ {\it anchored} if it runs through $\h_\d^\pm$,
                and otherwise we call the curve {\it floating}.

\bgnP{P:oned} A connected, reduced almost alternating link diagram with more than one 
              dealternator is the diagram depicted in Figure $\ref{F:2dalt}$.
              \endPR    
              Let $\a$ be one of the dealternators of the link diagram. Then $\a$ is adjacent 
              to four crossings. Since the crossing change at another dealternator $\b$ makes 
              the link diagram alternating and the diagram is reduced, each of the four 
              crossings is $\b$. 
              \endR 

\bgnF  {\iclg[scale=.9, bb=368 404 461 460]{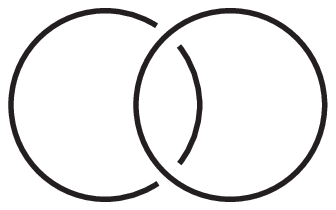}} \captn{} \label{F:2dalt} \endF

\bgnL{L:4ifc} Assume that $E$ is in basic position. 
              Let $C$ be an innermost curve of $E \cap S^\pm$ which is standard and floating.
              If  $C$ admits only trivial bridges, then $C$ has type $\Ga_1^\pm$.
              \endLR 
              We prove only the case when $C$ is positive, since the other case can be 
              shown similarly. Since $C$ is standard, $C$ meets a bubble or a segment.
              Moreover since $C$ is innermost, $C$ bounds a disk $D$ on $S^+$ whose interior 
              contains no positive curves. Thus the interior of $D$ does not contain the center 
              of an upper hemisphere, since otherwise $D$ contains the curve running through 
              the center of the upper hemisphere. 
              If $C$ meets bubbles in succession, then $D$ contains the center of one of the
              two upper hemispheres, since $C$ is floating. Thus $C$ meets a segment $\l$.
              Assume that $\l$ is $n$-adjacent to a bubble $\h$ and that $C$ runs through the 
              center of $\h^+$. Then $C$ meets a segment $\l'$ which is $p$-adjacent to $\h$, 
              i.e. $C$ has a neck, since $C$ is innermost and floating. Thus from Lemma 
              \ref{L:3neck}, $C$ has type $\Ga_1^+$, since $C$ admits only trivial bridges.
              Next assume that $C$ runs through $\l$. Then $\l$ is $p$-adjacent to a bubble $\h$
              which is not the dealternator, and $C$ runs through the center of $\h^+$,
              since $C$ is innermost and floating. Thus $C$ has type $\Ga_1^+$ as shown above.
              \endR 
 
\bgnL{L:4fcb} Assume that $E$ is in standard position. 
              If $E$ either has a floating curve or admits a non-trivial bridge,  
              then $E$ has a short arc of type ${\rm II}$ or ${\rm III}$.
              \endLR    
              We show only the case when $E$ has a positive floating curve or a non-trivial 
              positive bridge, since other cases can be shown similarly.
        \pvbn
              First consider the case when $E$ has a positive floating curve. 
              Take an innermost one $C$.  If $C$ admits only trivial bridges, then $C$ has 
              type $\Ga_1^+$ from Lemma \ref{L:4ifc}, and thus we are done. 
              If $C$ admits a non-trivial bridge $\g$, then operate the band move on $C$ 
              along $\g$ to have two positive curves $C'$ and $C''$.
              Here we have the following.

\bgnM{M:2ntb} We may take $\g$ so that $C'$ admits only trivial bridges and that 
              $C'$ and $C''$ are standard.
              \endMR
              Let $D$ be the subdisk of $E$ which $C$ spans in $B^+$.
              Then $C'$ spans a disk $D'=\Dgp \cup \t{D}'$, where $\t{D}'$ is 
              the component of $D-\g \times (-1,1)$ which contains $\g \times \{1\}$.
              Here we may assume that $D$ admits only trivial bridges in $\t{D}'$.
              Therefore if $C'$ admits a non-trivial bridge $\g'$, then $\g'$ intersects
              $\g \times \{1\}$ in one point $x$. Now let $\a_\g$ and $\b_\g$ be the ends 
              of $\g \times \{1\}$, and let $\a_{\g'}$ (resp. $\b_{\g'}$) be the end of 
              $\g'$ which is on the boundary of $\t{D}'$ (resp. of $\Dgp$).
              Let $\g_{\a_\g x}$ (resp. $\g_{x \a_{\g'}}$) be the subarc of 
               $\g \times \{1\}$ (resp. of $\g'$) whose ends are $x$ and $\a_\g$ 
              (resp. $\a_{\g'}$). Then $\g_{\a_\g x} \cup \g_{x \a_{\g'}}$ is
              a non-trivial bridge of $D$ in $\t{D}'$, which is a contradiction.
        \pvbn
              Since $\g$ is a bridge which is non-trivial, it is clear that
              $C'$ and $C''$ satisfy conditions (Et1-4). We can also see that
              $C'$ and $C''$ satisfy condition (Et5) by taking $\g$ so that
              $D$ admits only trivial bridges in $\t{D}'$ as above.
              \endR
        \pvbn
              From Claim \ref{M:2ntb} and Lemma \ref{L:4ifc}, $C'$ has type $\Ga_1^+$. 
              Therefore $C$ has a short arc of type II or IIIa, since $C$ is obtained  by 
              connecting $C'$ with $C''$ along a subarc of only one of the two short arcs. 
        \pvbn
              Next consider the case when $E$ has no floating curves but admits a 
              non-trivial bridge. Take a curve $C$ which admits a non-trivial bridge. 
              Then we can obtain from $C$ an innermost curve which is standard, floating, 
              and admits only trivial bridges by the band move along a non-trivial bridge. 
              Hence we know that $C$ has a short arc of type II or IIIa as above.
              \endR 

\bgnP{P:nosa} Assume that $E$ is in standard position. 
              Then $E$ is in special position if and only if $E$ has no short arcs.
              \endPR
              Note that if $E$ either has a floating curve or admits a non-trivial bridge, 
              then $E$ has a short arc from Lemma \ref{L:4fcb}. Assume that $E \cap S^\pm$ 
              has a curve $C$ which is not special. If $C$ has a neck, 
              then $C$ admits a non-trivial bridge or a short arc from Lemma \ref{L:3neck}.

\bgnM{M:nosa} If $C$ shares more than $1$ junction with an alternating segment, then
              $E$ admits a short arc.
              \endMR
              Assume that $C$ meets an alternating segment $\l$ in two junctions $a$ and $b$. 
              Here we consider the case when $C$ is positive. 
              The case when $C$ is negative can be shown similarly, and thus we omit it. 
              We see that a curve of $E \cap S^\pm$ can share at most $2$ junctions with 
              a segment, which are the ends of a boundary arc considering the orientation 
              of the curve and the segment, since the curve is a simple closed curve and 
              $E$ is in standard position. Thus $C$ runs through $\l$.
              Let $\h_1$ be the bubble which $\l$ is $n$-adjacent to. 
              Take the positive curve $C_1$ which runs through the center of 
              $\h_1^+$ and let $c$ be the junction of $C_1$ and $\l$. 
              Here we may assume that $b$ and $c$ are neighboring junctions. 
              Moreover then we may assume that $C$ neighbors $C_1$, i.e. there is no positive 
              curves on $S^+$ between $C$ and $C_1$, since otherwise $E$ has a floating curve. 
              Let $R$ be the region which contains the outside or isolated arc of $C$ with $b$  
              as an end of it. If $R$ has degree $2$, then obviously $C$ has a short arc in $R$. 
              If $R$ has degree no less than $3$, then take the bubble $\h_2$ of $R$ which 
              is $p$-adjacent to $\h_1$. Note that $\h_1$ does not have saddle-intersections, 
              since $C$ (resp. $C_1$) does not run through a side of $\h_1^+$ from condition
              (Et4) (resp. (Et3)) and $C$ neighbors $C_1$.
              Moreover we may assume that $\h_2$ is not at the dealternator, 
              since otherwise either $C$ or $C_1$ is floating. 
              Therefore $\l_{\h_1\h_2}$ is alternating, and thus has a junction. 
              Let $C_2$ be the positive curve with the closest junction to $\h_1$ on $\l_{\h_1\h_2}$. 
              Since $\h_1$ has no saddle-intersections and $C$ neighbors $C_1$, we have that 
              $C_2=C$ or that $C_2=C_1$.  
              In the former case $C_2$ admits a non-trivial bridge or has a short arc, and 
              in the latter case $C_2$ has a neck.
              \endR
        \pvbn
              Assume that $C$ is positive and that $C$ meets a non-alternating segment $\l$. 
              The case when $C$ is negative can be shown similarly, and thus we omit the proof.
              If $\l$ is $p$-adjacent to the dealternator, then $C$ runs through $\l$.
              Then $C$ does not run through a side of $\h_\d^+$ from condition (Et4).
              Thus $C$ either is floating or has a neck.
              Next if $\l$ is $n$-adjacent to the dealternator, then there is a negative 
              curve which runs through $\l$, and thus we can be done similarly.      
        \pvbn
              Conversely assume that $E$ has a short arc. Then take a crossing $x$ which 
              admits a short arc and take the closest short arc $\q$ to $x$. 
              Let $C_\q^\pm$ be the positive/negative curve which contains $\q$ and 
              let $\l$ be the end segment of $\q$ which is $p$-adjacent to $x$. If $\l$ is 
              non-alternating, then neither $C_\q^+$ nor $C_\q^-$ satisfies condition (Ep3). 
              Next assume that $\l$ is alternating. 
              If $\q$ has type I or IIIa, then $C_\q^+$ has two junctions with $\l$, since 
              $C_\q^+$ does not run through the center of $\h_x^-$. 
              Thus $C_\q^+$ does not satisfies condition (Ep2). If $\q$ has 
              type II or IIIb, then $C_\q^-$ does not satisfies condition (Ep1) or (Ep2).  
              \endR


\section{A spanning disk of the trivial knot in almost alternating position I}  
\label{S:CUT}                                                                   

               Let $K$ be the trivial knot in strongly reduced almost alternating position.
               Let $E$ be a spanning disk for $K$ in basic position with minimal complexity.
              Thus $K$ has only one dealternator from Proposition \ref{P:oned}, and 
                   $E$ is in standard position from Proposition \ref{P:MC/ST}.
               Moreover $K$ is prime from the following proposition. 
 
\bgnP{P:ctpr} A connected, reduced almost alternating diagram of a trivial link is prime.
              \endPR    
              Let $L$ be a non-prime, connected, reduced almost alternating link diagram.
              From Proposition \ref{P:oned}, we have that $L$ has only one dealternator.
              Thus $L$ can be decomposed into a connected alternating diagram $L'$ and a 
              connected almost alternating diagram $L''$ such that $L'$ is reduced. 
              Then $L'$ does not represent a trivial link from Theorem \ref{T:AL} or from
              Theorem \ref{T:AK}. This implies that $L$ does not represent a trivial link
              (see \cite{BZ-85} Corollary 7.5 (b)).
              \endR     
                        
\bgnL{L:5dsc}  The dealternator of $K$ does not admit a short cut.
               \endLR 
               Assume that the dealternator admits a short cut $\m$. Then the dealternator
               of $K^\m$ is a nugatory crossing. Since $K$ is prime, $K^\m$ is connected.
               Thus we obtain a connected alternating diagram of the trivial $2$-component 
               link from $K^\m$ by the Reidemeister move of type I.
               This contradicts to Theorem \ref{T:AL}.                   
               \endR   

\bgnP{P:para}  Non-equivalent disjoint short cuts of $E$ are parallel.
               \endPR
               Assume that $E$ admits non-equivalent disjoint short cuts $\m$ and $\m'$. 
               Let $\wE_1$, $\wE_2$ and $\wE_3$ be subdisks of $E$ such that
               $\wE_1 \cup \wE_2 \cup \wE_3 =  E $, $\wE_1 \cap \wE_2 = \m $  and 
                          $\wE_2 \cap \wE_3 = \m'$. Then we obtain disks 
               $E_1$, $E_2$ and $E_3$ by cut surgeries along $\m$ and $\m'$ such that: 
\pvan          $E_1 =(\wE_1 - \m  \times (-1,0) ) \cup \Dmn$;     \\%
               $E_2 =(\wE_2 - \m  \times ( 0,1) - \m' \times (-1,0) ) \cup \Dmp \cup \Ddn$; and \\%
               $E_3 =(\wE_3 - \m' \times ( 0,1) ) \cup \Ddp$.     
         \pvan 
               Let $K_i=$ $\pE_i$ ($i=1,2,3$). Note that $K^\m$ and $K^{\m'}$ are connected, 
               since $K$ is prime. Thus $(K^\m)^{\m'}$ is a disconnected almost altnernating 
               diagram of the trivial $3$-component link consisting of two connected components 
               of an almost alternating diagram $K_1 \cup K_3$, and a trivial or alternating 
               diagram $K_2$ 
               from Lemma \ref{L:5dsc}, Theorem \ref{T:AAL} (1) and Theorem \ref{T:AL}.                
               Since $\m$ and $\m'$ do not have the ends on same segments, $K_2$ has 
               a crossing and thus $K_2$ is a coiled diagram from Theorem \ref{T:AK}. 
         \pvcn 
               Let $x_1$, $\cdots$, $x_k$ be the crossings of $K_2$, where $x_i$ and $x_{i+1}$
               belong to a common region $R_i$ of degree $2$, and $x_1$ and $x_k$ admits $\m$ 
               and $\m'$ in $E$, respectively. 
               We claim here that each $\h_{x_i}$ has no saddle-intersections. 
               From Lemma \ref{L:3sab}, we know that neither $\h_{x_1}$ nor $\h_{x_k}$ has
               saddle-intersections. Thus assume that $k$ is no less than $3$ and 
               that $\h_{x_2}$ has a saddle-intersection. Then there is a positive curve 
               which runs through a side of $\h_{x_2}^+$ and goes into $R_1$.
               Then the curve runs through a side of $\h_{x_1}^+$, since $R_1$ has 
               degree $2$ and $E$ is in standard position. However this contradicts 
               that $\h_{x_1}$ has no saddle-intersections. Now the claim holds by an induction.  
               Therefore the curves of $(E_1 \cup E_3) \cap S^\pm$ are away from $K_2$. 
               Then $E_2$ is coiled, i.e. each curve of $E_2 \cap S^\pm$ has type 
               $\Ga_1^\pm$ or $\Ga_2^\pm$, since $E$ has a minimal complexity. 
               Hence $\m$ and $\m'$ are parallel.   
               \endR
 
\bgnC{C:5csc}  A curve of $E \cap S^\pm$ does not admit non-equivalent disjoint short cuts.
               \endC \End

\bgnC{C:5fcb}  An innermost floating curve $C$ of $E\cap S^\pm$ has 
               type $\Ga_1^\pm$ or $\Ga_2^\pm$.
               \endCR
               If $C$ admits only trivial bridges, then $C$ has type $\Ga_1^\pm$ from 
               Lemma \ref{L:4ifc}. If $C$ admits a non-trivial bridge, then sufficiently many 
               band moves on bridges split $C$ into a set of curves each of which has type 
               $\Ga_1^\pm$. Since $C$ does not admit non-equivalent disjoint short cuts from 
               Corollary \ref{C:5csc}, the set consists of two curves and the bridge of 
               $C$ is a short cut. Therefore $C$ has type $\Ga_2^\pm$, since $K$ is
               strongly reduced.
               \endR            

\bgnL{L:5xsa}  A crossing of $K$ does not admit two short arcs in a same region.
               \endLR
               Assume that a crossing admits a short arc $\q$ of type I or II. 
               Then the two surgered segments of $K^\q$ belong to different components 
               of the trivial $2$-component link, implying the conclusion.    
         \pvcn                     
               Next assume that a crossing admits two short arcs $\q$ and $\q'$ in a region. 
               Here we assume that there are no short arcs on the region between $\q$ and $\q'$.
               Thus both of $\q$ and $\q'$ have the same type of IIIa or IIIb from the above.
               We consider only the former case, since the latter case can be shown similarly.
         \pvcn    
               Let $x$ be a crossing which admits short arcs $\q=a_1a_2$ and $\q'=b_1b_2$ 
               of type IIIa in a region $R$, where $a_i$ and $b_i$ are junctions on a segment 
               $\l_i$ ($i=1,2$) and $\q'$ is closer to $x$ than $\q$ on $R$. Let $\g=c_1c_2$ 
               be a positive trivial bridge $\q' \times (-1)$ at $x$ such that $c_i$ 
               is on $C_{\q'} \cap \l_i$, where $C_{\q'}$ is the positive curve containing $\q'$. 
               Let $D$ be the subdisk of $E$ bounded by $\q$ and $\g$. Take a disk $D'$ 
               in $B^+$ such that $\pD = \g \cup \l_{c_1a_1} \cup \q \cup \l_{a_2c_2}$ and
               $D' \cap E = D' \cap D$, where $\l_{c_ia_i}$ (resp. $\l_{a_ic_i}$) is 
               the subsegment of $\l_i$ with $a_i$ and $c_i$ as its ends.
               Then replace $D$ with $D'$ to have 
               another spanning disk for $K$, which is clearly in basic position with a 
               fewer complexity than that of $E$. This contradicts the minimality of $E$.      
               \endR  

\bgnL{L:5dsd}  The bubble at the dealternator of $K$ has a saddle-intersection.
               \endLR  
               Assume otherwise and let $C_\d^\pm$ be the curve which runs through the center 
               of $\h_\d^\pm$. Then $C_\d^+$ and $C_\d^-$ are the only curves which run 
               through $\h_\d$. Let $D_i^\pm$ be a disk bounded by $C_\d^\pm$ on $S^\pm$ 
               such that $D_1^\pm \cap D_2^\pm =C_\d^\pm$ and
                         $D_1^\pm \cup D_2^\pm =S^\pm$ ($i=1,2$).
               We claim here that the interior of $D_1^\pm$ or of $D_2^\pm$ contains no 
               positive/negative curves. Assume otherwise and let $C_i^\pm$ be an innermost curve 
               in the interior of $D_i^\pm$. Since any curve other than $C_\d^\pm$
               is floating, $C_i^\pm$ admits a short cut $\g_i^\pm$ from Corollary \ref{C:5fcb}. 
               However then $\g_1^\pm$ 
               and $\g_2^\pm$ are not parallel, since we have $C_\d^\pm$ on $S^\pm$ between 
               $\g_1^\pm$ and $\g_2^\pm$. This contradicts Proposition \ref{P:para}.
         \pvbn 
               Therefore $C_\d^\pm$ bounds a disk $D^\pm$ on $S^\pm$ whose interior contains 
               no positive/negative curves. Let $x$ be the crossing which is $p$-adjacent to 
               the dealternator so that $\l_{x\d}$ meets $D^+$. If $C_\d^+$ meets $\l_{x\d}$, 
               then the dealternator admits a short cut from Lemma \ref{L:3neck}. 
               This contradicts Lemma \ref{L:5dsc}. Thus $\l_{x\d}$ is contained in the interior 
               of $D^+$. Since the interior of $D^+$ contains no positive curves, $\l_{x\d}$ has
               no junctions and then $C_\d^-$ runs through the center of $\h_x^-$. Thus $\h_x$ 
               has no saddle-intersections, since $C_\d^-$ bounds $D^-$. Hence $C_\d^+$ runs 
               through the center of $\h_x^+$. Then we can take a disk $\De$ in $B^+$ such that 
               $\prtl \De=\a \cup \b \cup \l_{x\d} \cup \c$, where $\De \cap E=\a$, 
               $\De \cap S^+=     \b \cup \l_{x\d} \cup \c$ and $\b$ (resp. $\c$) is on $\h_x^+$ 
               (resp. $\h_\d^+$). Thus both of $\g_1=\a \times \{-1\}$ and $\g_2=\a \times \{1\}$ 
               are positive bridges. Then $(K^{\g_1})^{\g_2}$ is a connected almost alternating 
               diagram of the trivial $3$-component link, since $K$ is strongly reduced.
               However this contradicts Theorem \ref{T:AAL} (1).    
               \endR

\bgnL{L:5dsa}  The dealternator of $K$ does not admit a short arc.
               \endLR         
               Let $x$ be a crossing which admits a short arc. 
               Then $\h_x$ does not have a saddle-intersection from Lemma \ref{L:3sab}.
               Thus $x$ is not the dealternator from Lemma \ref{L:5dsd}.
               \endR


               \subsection{A spanning disk with a short arc} 
\label{SS:c1}
\bgnP{P:cut1}  If $E$ has a short arc $\q$ of type ${\rm I}$, then $K^\q$ is the trivial 
               $2$-component link in prime, strongly reduced almost alternating position.
               \endPR  
               Let $x$ be a crossing which admits $\q$ and take a look at the diagram for $K^\q$
               in Figure \ref{F:sarc}. Since $x$ is not the dealternator from Lemma \ref{L:5dsa},
               we can see that $K^\q$ is almost alternating and strongly reduced. Moreover 
               since $K^\q$ is clearly connected, $K^\q$ is prime from Proposition \ref{P:ctpr}.
               \endR

\bgnL{L:5sa1}  If $E$ has a short arc of type ${\rm I}$, then $E$ has no other short arcs.
               \endLR
               If $E$ has a short arc $\q$ of type I and another short arc $\q'$, then 
               $(K^\q)^{\q'}$ is a connected almost alternating diagram of the trivial 
               $3$-component link from Lemma \ref{L:5dsa} and Proposition \ref{P:cut1}. 
               This contradicts Theorem \ref{T:AAL} (1). 
               \endR

\bgnP{P:a2sc}  $E$ has a short arc of type $\rm{II}$ or $\rm{III}$ if and only if 
               $E$ admits non-equivalent disjoint short cuts.
               \endPR
               If $E$ admits a short cut, then $E$ has a short arc of type II or III from 
               Lemma \ref{L:4fcb}, since a short cut is a short arc of type III or a short bridge.
         \pvbn
               Next assume that $E$ has a crossing $x$ which admits a short arc $\q$ of type II.
               Let $C_\q^\pm$ be the positive/negative curve containing $\q$ and let 
               $\l_i^+$ (resp. $\l_i^-$) be the segment which is $p$-adjacent (resp. $n$-adjacent)  
              to $x$ ($i=1,2$), where $\l_1^+$ and $\l_1^-$ are the end segments of $\q$. 
               Since $C_\q^\pm$ runs through the center of $\h_x^\pm$ from Lemma \ref{L:5xsa},             
               $C_\q^\pm$ has a short cut $\m^\pm$ whose ends are on $\l_1^\pm$ and $\l_2^\mp$.
               Then  $\m^+$ and $\m^-$ are non-equivalent and disjoint, 
               since $\m^+$ and $\m^-$ have different end segments.
         \pvbn
               If $E$ has a floating curve, then an innermost floating curve has type $\Ga_1$ 
               or $\Ga_2$ from Corollary \ref{C:5fcb}, and thus a short arc of type II. Therefore 
               we are done from the above. We complete the proof by showing the following claim.

\bgnM {M:a2bc} If $E$ has a short arc of type {\rm III}, then $E$ has a short arc of type {\rm II}.
               \endMR
               We show only the case when $E$ has a short arc of type IIIa, since the other case 
               can be shown similarly. In addition we may assume that $E$ has no floating curves
               from the above. Let $x$ be a crossing which admits a short arc $\q$ of type IIIa  
               and let $C_\q$ be the positive curve containing $\q$. 
               Then $C_\q$ runs through the center of $\h_x^+$ from Lemma \ref{L:5xsa}.
               Let $y$ and $z$ be the crossings which are $p$-adjacent to $x$ 
               with segment $\l_{yx}$ containing an end of $\q$.  
         \pvbn
               Assume that $y$ is the dealternator. Then $C_\q$ runs through the center of 
               $\h_y^+$ from condition (Et4). However then, $C_\q$ has a neck and thus has
               a short cut from Lemma \ref{L:3neck}, contradicting Lemma \ref{L:5dsc}.
         \pvbn
               Assume that $z$ is the dealternator. 
               Consider the case when $C_\q$ runs through the center of $\h_z^+$. 
               From Lemma \ref{L:5dsd}, $\h_z$ has a saddle-intersection. 
               Then a positive curve running through $\h_z^+$ on the side of 
               $\l_{zx}$ runs through $\h_x^+$ on the side of $\l_{zx}$.   
               However this contradicts Lemma \ref{L:3sab}. 
               Next consider the case when $C_\q$ runs through a side of $\h_z^+$. 
               If there is a positive curve running through a side of $\h_z^+$ closer to $\l_{zx}$ 
               than $C_\q$, then we obtain a contradiction to Lemma \ref{L:3sab} as above. 
               Otherwise $C_\q$ bounds a disc $D$ on $S^+$ such that the center of $\h_z^+$ 
               is in the interior of $D$ and $S^+-D$ has no positive anchored curves.
               Then the curve running through the center of $\h_y^+$ is in $S^+ -D$,
               and thus floating. This contradicts the assumption.           
         \pvbn
               Now assume that neither $y$ nor $z$ is the dealternator. Thus both of $\l_{yx}$ and 
               $\l_{zx}$ are alternating segments. Let $a_1$, $a_2$ and $a_3$ be consecutive 
               junctions on $\l_{yx}$ such that $\l_{a_1a_2}$ is $\l_{yx} \cap C_\q$ with $a_2$
               an end of $\q$. Then the positive curve $C_a$ which runs through $a_3$
               is not $C_\q$ but neighbors $C_\q$, i.e. there are no positive curves on $S^+$
               between $C_a$ and $C_\q$, since $E$ has no floating curves.
               Next let $C_b$ be the positive curve which runs through the closest junction $b$ 
               to $\h_x$ on $\l_{zx}$. Then $C_b$ is not $C_\q$, since otherwise $C_\q=C_b$ admits
               non-equivalent disjoint short cuts at $x$, contradicting Corollary \ref{C:5csc}.
               Thus $C_b$ neighbors $C_\q$, since no curves run through a side of $\h_x$ and
               $E$ has no floating curves. Therefore $C_\q$ 
               has two neighbors in a same component of $S-C_\q$, and thus we have that $C_a=C_b$.
               However this is impossible, since the segment running through $\h_y^+$ and
               the segment running through $\h_z^+$ belong to different components in $K^\q$.
               \endR
               \endR  

  
               \subsection{A spanning disk with non-equivalent disjoint short cuts}
\label{SS:c2}
               Let $X$ be a set of mutually non-equivalent disjoint short cuts such that 
               any short cut of $E$ either intersects or is equivalent to an element of $X$. 
               From Proposition \ref{P:para}, there is a pair of short cuts, say $\g_l$ and $\g_r$,
               which bounds a subdisk $\Om$ of $E$ containing all the elements of $X$.
               We define the {\it extract surgery on} $E$ as the operation of getting rid 
               of $(\Om^{\g_l})^{\g_r}$ from $(E^{\g_l})^{\g_r}$, and denote the result
               by $E^*$ and $\pE^*$ by $K^*$.
\bgnP{P:exsg}  The extract surgery on $E$ is well-defined.
               \endPR
               Assume that there is a short cut $\g'$ which is not equivalent to any element of $X$.
               Then there is a crossing $x$ which admits a short cut $\g$ of $X$ which intersects $\g'$. 
               Let $C$ be the curve of $E \cap S^\pm$ admitting both of $\g$ and $\g'$. 
               If $\g$ is neither $\g_l$ nor $\g_r$, then $C$ has type $\Ga_1^\pm$ or $\Ga_2^\pm$
               from Proposition \ref{P:para}. However this is a contradiction, since a curve
               with type $\Ga_1^\pm$ or $\Ga_2^\pm$ does not admit non-equivalent short cuts.
               Thus assume that $\g$ is $\g_l$. Then $C$ has a neck and a short arc of type II
               at $x$ from Proposition \ref{P:para}. Thus $\g'$ has a common end segment with $\g$.
               Let $a$ and $b$ (resp. $a'$ and $b'$) be the ends of $\g$ (resp. $\g'$), where
               $a$ and $a'$ are on a same segment. Then we can obtain non-equivalent disjoint 
               short cuts, one of which has ends $a$ and $b'$ and the other has ends $a'$ and $b$ 
               by smoothing the intersection of $\g$ and $\g'$.
               However this contradicts Corollary \ref{C:5csc}.
               \endR
\nidt
               Take an arc $\t{\psi}$ on $\Om$ which connects $\g_l$ and $\g_r$.
               Let $\psi$ be a projection of $\t{\psi}$ on $S^+ \cap S^-$ and call $\psi$
               a {\it band-trace} for $E$.
\bgnP{P:g2hl}  Let $L$ be a connected, reduced almost alternating diagram of the trivial 
               $2$-component link. If $L$ is not strongly reduced, then $L$ is the diagram 
               of Figure $\ref{F:2dalt}$.
               \endPR
               Apply the Reidemeister move of type II to $L$ to have another diagram $L'$.
               Then $L'$ is alternating or trivial. In the former case $L'$ is disconnected
               from Theorem \ref{T:AL}. Since $L$ is prime from Proposition \ref{P:ctpr},
               each component of $L'$ is reduced. This contradicts Theorem \ref{T:AK}.
               In the latter case we have the conclusion.
               \endR

\bgnP{P:cut2}  Assume that $E$ admits non-equivalent disjoint short cuts. Then $K^*$ is 
               the trivial $2$-component link in prime, strongly reduced almost alternating 
               position and $E^*$ is in special position.
               \endPR
         \nidt
               Take a look at $E$ and recall notations in the definition of the extract surgery.  
               Since $K$ is prime, there uniquely exists a region which contains the two end
               segments of $\g_i$ $(i=l,r)$. We denote the region by $R_i$.
               Define crossings $x_i$, 
               $y_i$ and $z_i$ of $(E-\Om)\cap S$ if ${\rm deg}R_i=2$, and of $\pR_i$ if 
               ${\rm deg}R_i \geq 3$ as follows: $x_i$ be the crossing which admits $\g_i$; and 
               $y_i$ (resp. $z_i$) be the crossing which is $p$-adjacent (resp. $n$-adjacent) to 
               $x_i$. Let $C_{\g_i}$ be the curve which admits $\g_i$, where we take the one 
               which is not contained in $\Om$ if $\g_i$ is a short arc.
         \pvcn
                Since $K$ is almost alternating and prime, we have that $x_r \neq y_l$, $z_l$ 
                and equivalently that $x_l \neq y_r$, $z_r$. Next assume that $y_l=z_r=\d$ or 
                that $y_r=z_l=\d$. Take a look at the region $R$ whose boundary crossings are $\d$,
                $x_l$, $x_r$ and other crossings of $\Om$. From Lemma \ref{L:5dsd}, there is a
                positive curve $C$ which runs through a side of $\h_\d^+$ and goes into $R$.              
                However then $C$
                meets neither $\l_{\d x_l}$ nor $\l_{\d x_r}$ from condition (Et4), and $C$
                does not meet $\Om$ from the proof of Proposition \ref{P:para}.
                Hence the dealternator is not adjacent to both of $x_l$ and $x_r$, 
                    and thus we may assume that $x_r$ is not adjacent to the dealternator. 
                In addition, we may assume that $C_{\g_l}$ is positive, since the other case 
                can be shown similarly. 
         \pvcn
                Here we define $C_l^+$ and $C_l^-$. Let $C_l^+$ be $C_{\g_l}$.
                If $R_l$ has degree $2$ (resp. has degree no less than $3$ and $\l_{x_l z_l}$ 
                has a junction), then let $C_l^-$ be the curve of $(E-\Om)\cap S^-$ 
                sharing a junction with $C_l^+$ on the segment facing $R_l$ which is $p$-adjacent 
                (resp. $n$-adjacent) to $x_l$; and if $R_l$ has degree no less than $3$ and
                $\l_{x_l z_l}$ has no junctions, then let $C_l^-$ be the curve which runs 
                through $\h_{z_l}$ the closest side to $\l_{x_l z_l}$, where such a curve
                exists from Lemma \ref{L:5dsd}, since now $z_l$ is the dealternator.
         \pvcn
                Next define $C_r^+$ and $C_r^-$ as follows; if $C_{\g_r}$ is a positive curve, then 
                let $C_r^+$ be $C_{\g_r}$ and let $C_r^-$ be the curve of $(E-\Om)\cap S^-$ sharing 
                a junction with $C_r^+$ on the segment facing $R_r$ which is $n$-adjacent (resp. 
                $p$-adjacent) to $x_r$ if $R_r$ has degree no less than $3$ (resp. degree $2$);
                and if $C_{\g_r}$ is a negative curve, then let $C_r^-$ be $C_{\g_r}$ and let 
                $C_r^+$ be the curve of $(E-\Om)\cap S^+$ sharing a junction with $C_r^-$ on the 
                segment facing $R_r$ which is $p$-adjacent (resp. $n$-adjacent) to $x_r$ if $R_r$ 
                has degree no less than $3$ (resp. degree $2$).

\bgnM  {M:c2fc} $E^*$ has no floating curves.
                \endMR
                It is sufficient to show that $E$ does not have a second-innermost floating
                curve, since every innermost floating curve of $E$ belongs to $\Om \cap S^\pm$. 
                If $E \cap S^\pm$ has no innermost floating curves, then we are done. Thus 
                assume otherwise. From the proof of Proposition \ref{P:para}, no curves of 
                $(\overline{E-\Om}) \cap S^\pm$ run through the bubble at a crossing 
                of $\Om \cap S$. Therefore if $E$ has a second-innermost floating curve, 
                then it is $C_l^\pm=C_r^\pm$. However this is impossible, since: 
                $C_l^+ \cap (E-\Om)$ and $C_r^+ \cap (E-\Om)$ belong to different components
                of $E-\Om$; $C_l^- \cap (E-\Om)$ and $C_r^- \cap (E-\Om)$ belong to different 
                components of $E-\Om$ if $R_l$ has degree $2$ or if $R_l$ has degree no less
                than $3$ and $\l_{x_l y_l}$ has a junction; and $C_l^-=C_r^-$ runs through
                the dealternator if $R_l$ has degree no less than $3$ and $\l_{x_l y_l}$ has 
                no junctions.
                \endR
         \pvan
                Here assume that $K^*$ is reduced and $E^*$ is in standard position. 
                Consider $E^*$ with a band-trace $\psi$. Here note that $K^*$ is a 
                connected almost alternating diagram of the trivial $2$-component link. 
                Assume that $K^*$ is not strongly reduced. Then $K^*$ is the diagram  
                of Figure \ref{F:2dalt} from Proposition \ref{P:g2hl}.
                Since each of the four regions of $S^+$ with $K^*$ is a trivial clasp
                and $\psi$ is in one of the four regions, $K$ is not strongly reduced, either.
                Thus reducedness of $K^*$ implies strongly reducedness of $K^*$. 
                Next if $K^*$ has a short arc $\q$, then
                $E^*$ has a short cut $\g$ from the assumption, Lemma \ref{L:5sa1} and
                Proposition \ref{P:a2sc}, since the extract surgery does not creat
                new non-boundary arcs for $E^*$. Therefore $(K^*)^\g$ is a connected almost
                alternating diagram of the trivial $3$-component link from Lemma \ref{L:5dsc},
                since $K^*$ is prime from Proposition \ref{P:ctpr}.
                This contradicts Theoerm \ref{T:AAL} (1). 
                Thus $E^*$ is in special position from Proposition \ref{P:nosa}. 
                Therefore it is sufficient to show the following two claims.

\bgnM  {M:c2rd} $K^*$ is reduced.
                \endMR
                It is sufficient to show that each of $R_l$ and $R_r$ has degree no less 
                than $3$. Take a look at $E \cap S^\pm$. Since the proof is similar to the 
                proof of Claim \ref{M:a2bc}, we omit the detail in the following.
         \pvbn
                Assume that $R_l$ has degree $2$. Then $y_l$ is not the dealternator, since 
                otherwise $C_r^+$ is floating, implying a contradiction to Claim \ref{M:c2fc}. 
                Thus $\l_{x_l y_l}$ is alternating, and we have a positive curve which has the 
                closest junction to $x_l$ on segment $\l_{x_ly_l}$. Then we obtain 
                a contradiction by considering the curve with $C_l^+$ and $C_r^+$.
         \pvbn
                Assume that $R_r$ has degree $2$ and $C_{\g_r}$ is positive. Then we have a 
                positive curve which has the closest junction to $x_r$ on segment $\l_{x_ry_r}$, 
                since $y_r$ is not the dealternator. Then we obtain a contradiction by 
                considering the curve with $C_l^+$ and $C_r^+$.
         \pvbn
                Now assume that $R_l$ has degree no less than $3$, and that $R_r$ has degree $2$ 
                and $C_{\g_r}$ is negative.  
                Note that $C_r^+$ runs through the center of $\h_{x_r}^+$. 
                Let $C$ be the negative curve which shares a junction with $C_r^+$ on 
                segment $\l_{x_rz_r}$. Assume that $\l_{x_lz_l}$ has no junctions. 
                Since no curves of $(E-\Om) \cap S^+$ run through a crossing of $\Om \cap S$, 
                $C_r^-=C_l^-$ and thus
                $C_r^-$ runs through $\h_{z_l}^-=\h_\d^-$ on the closest side to $\l_{z_lx_l}$. 
                Thus $C$ is floating, implying a contradiction to Claim \ref{M:c2fc}. 
                If $\l_{x_lz_l}$ has a junction, then we obtain a 
                contradiction as above by considering the curves $C_l^-$, $C_r^-$ and $C$.
                \endR
\bgnM{M:c2st}   $E^*$ is in standard position.
                \endMR
                Take a look at $E \cap S^\pm$.
                First we show that $y_l$ is not the dealternator. Assume otherwise. 
                Note that $C_l^+$ is anchored from Claim \ref{M:c2fc}. 
                If $C_l^+$ runs through a side of $\h_{y_l}^+$, then $C_l^+$ does not satisfies 
                condition (Et4). If $C_l^+$ runs through the center of $\h_{y_l}^+$, then $C_l^+$ 
                admits a neck at $y_l$, and thus $y_l$ admits a short cut from 
                Lemma \ref{L:3neck}. This contradicts Lemma \ref{L:5dsc}.
         \pvbn  
                Since $C_{\g_l}$ and $C_{\g_r}$ are the only curves which are changed by
                the extract surgery, it is sufficient to show that these two curves are standard
                after the surgery. This is done by showing that $C_{\g_l}$ runs 
                through the center of $\h_{y_l}^+$ in $E$ and that $C_{\g_r}$
                runs through the center of $\h_{y_r}^+$ (resp. $\h_{z_r}^-$) if $C_{\g_r}$ 
                is positive (resp. negative). If any of these claims does not hold, we can
                show that $y_l$, $y_r$ or $z_r$ admits a short cut from Lemma \ref{L:5sa1}
                by following the proof of Claim \ref{M:nosa}. However this is a contradiction.
                \endR
                \endR


\section{A spanning disk of the trivial knot in almost alternating position II} 
\label{S:GLU}

               Let $K$ be the trivial knot in strongly reduced almost alternating position
               and let $E$ be a spanning disk for $K$ in basic position with minimal complexity.
               We assume here that $E$ has a short arc $\q$. Then after a proper surgery, 
               $K^*$ is the trivial $2$-component link in prime, 
               strongly reduced almost alternating position
               from results in Section \ref{S:CUT}, where we denote $K^\q$ by $K^*$ if $\q$ 
               has type I, since it causes no contradiction from Lemma \ref{L:5sa1}. Therefore 
               $K^*$ has a flyped tongue from Theorem \ref{T:AAL} (2). We show in this section
               that we can operate the inverse of the surgery on $E^*$ without harming the flyped 
               tongue of $K^*$.  
         \pvcn
               Define the {\it left} and {\it right} side of a non-alternating segment $\ldq$ by 
               running along $\ldq$ from crossing $q$ to the dealternator $\d$. Denote the region 
               which faces $\ldq$ from the left (resp. right) side by $O^q_l$ (resp. $O^q_r$). 
               Denote the region sharing with $O^q_i$ a segment ($\neq \ldq$) which is adjacent 
               to $\d$ (resp. $q$) by $P^q_i$ (resp. $Q^q_i$) ($i=l$, $r$). 
               We say that $\ldq$ has a {\it flype-component} on the left (resp. right) side or 
               a flype-component with a {\it flype-crossing} $x$ if $P^q_l$ 
               and $Q^q_l$ (resp. $P^q_r$ and $Q^q_r$) share a common crossing $x$.  
               If $\ldq$ has flype-components on both sides, $K^*$ has a flyped tongue. 
               Then we call the pair of $O^q_l$ and $O^q_r$ the {\it core} of a flyped tongue. 
         \pvcn
               Take a look at a flype-component of $K^*$ with a flype-crossing $x$. 
               We may denote $O^q_i$, $P^q_i$ and $Q^q_i$ by $O^q_x$, $P^q_x$ and $Q^q_x$ 
               ($i=$ $l$ or $r$). We call the $2$-tangle $W^q_x$ with $\d$, $q$ and $x$ 
               as its ends, the {\it flype-tangle of} $\ldq$ with $x$. 
               We say that $W^q_x$ is {\it trivial} if 
               $W^q_x$ consists of two segments $\l_{x\d}$ and $\l_{xq}$. 
               In the following we omit $q$ of regions unless we need to emphasize.


\pvcn
               \subsection{A spanning disk with a short arc of type I}
\label{SS:g1}

\bgnP{P:glu1}  If $E$ has a short arc $\q$ of type $\rm{I}$, then $K$ has a flyped tongue.
               \endPR
               From Proposition \ref{P:cut1}, 
               $K^*$ is the trivial $2$-component link in prime, 
               strongly reduced almost alternating position. Thus $K^*$ has a 
               non-alternating segment $\ldq$ which has a flype-component with crossing 
               $x_l$ (resp. $x_r$) on the left (resp. right) side from Theorem \ref{T:AAL} (2). 
               Consider the inverse operation of the cut surgery along $\q$
               paying attention only on $K^*$, which can be 
               regarded as an operation of smoothing one of the two crossings of a non-trivial
               clasp $\Si$ of $K^*$. If none of $\d$, $q$, $x_l$ and $x_r$ belongs to $\Si$, 
               then we see that $K$ also admits the flype-components and thus we are done. 
         \pvcn
               Since the dealternator cannot belong to $\Si$, it is sufficient to consider the 
               cases when $q$ or $x_l$ belongs to $\Si$. Assume that $q$ belongs to $\Si$. 
               Since $K^*$ is strongly reduced, each region 
               facing $\ldq$ has degree no less than $3$. Thus we may assume that $Q_l$ is $\Si$. 
               Then flype-tangle $W_{x_l}$ is trivial, and thus $K$ is not strongly
               reduced no matter which crossing of $q$ and $x_l$ we smooth.
         \pvcn
               Next assume that $x_l$ belongs to $\Si$ but $q$ does not. Since $\d$ does not
               belong to $\Si$, one of the two regions ($\neq P_l$, $Q_l$) which has $x$ is $\Si$.
               In either case, $P_l$ and $Q_l$ share the other crossing $y$ of $\Si$. 
               Then $y$ is another flype-crossing of $\ldq$ on the left side. Therefore $K$ 
               has a flype-component of $\ldq$ with $x_l$ if we smooth $y$, or with $y$ if we 
               smooth $x_l$. 
               \endR

\pvcn
               \subsection{A spanning disk with non-equivalent disjoint short cuts}
\label{SS:g2}
               Next we consider the case when $E$ has a short arc of type II or III, 
               i.e. admits non-equivalent disjoint short cuts from Proposition \ref{P:a2sc}. 
               Then $K^*$ is the trivial $2$-component link in prime, 
               strongly reduced almost alternating position
               and $E^*$ is in special position from Proposition \ref{P:cut2}. 
               We consider only the case when $K^*$ has a flyped tongue as the left of 
               Figure \ref{F:6lgs}, since the other case can be shown similarly.
               Since $E^*$ is in special position, $E^*$ admits neither a floating curve nor 
               a non-trivial bridge from Lemma \ref{L:4fcb} and Proposition \ref{P:nosa}.
               Thus first we have the following.

\bgnL{L:g2cr}  A curve of $E^* \cap S^\pm$ meets a region in an arc.          \endLR
               Otherwise $E^*$ admits a non-trivial bridge.                   \endR
         \pvcn
               Since $E^*$ does not have a floating curve, curves with a same sign are concentric
               on $S^\pm$. Let $C_\d$, $C_1$ and $C_2$ be curves of $E \cap S^\pm$, where $C_\d$ 
               is the curve running through the center of $\h_\d^\pm$. Then we say that 
               $C_\d > C_1 > C_2$ if $C_1$ bounds a disk in a component of $S^\pm-C_\d$ 
               which contains $C_2$. If $C_1$ and $C_2$ are in different components, 
               we say that $C_1>C_\d>C_2$ or $C_2>C_\d>C_1$.
         \pvcn
               In the following we denote by $C_x$ the positive curve which runs through the 
               center of $\h_x^+$ of a crossing $x$. From Lemma \ref{L:g2cr}, we know a curve 
               precisely if we are given which crossings and how the curve runs through, since
               $K^*$ is prime and $E^*$ is in special position.
               Thus we may denote a curve only by giving crossings with order which the curve runs 
               through, where we denote a crossing by itself (resp. itself with a bar on top) if
               the curve runs through a side (resp. the center) of the upper hemisphere of the crossing.
               We denote an arc of a positive curve similarly, e.g. an inside arc by $\c_{xy}$;
               an outside arc by $\c_{\bar{x}y}$; and an isolated arc by $\c_{\bar{x}\bar{y}}$.

\bgnL{L:g2cc}  Let $R$ be a region which has the dealternator $\d$. 
               Let $y_1$ be the crossing of $R$ which is $p$-adjacent to $\d$ and 
               let $y_i$ be the crossing of $R$ which is $n$-adjacent to $y_{i-1}$ 
               $(i=2, 3,\cdots,n)$ so that $y_n$ is $n$-adjacent to $\d$. Then we have that 
               $C_{y_j}$ contains an outside arc $\c_{\bar{y}_j\d}$ in $R$ and that
               $C_{y_n} =C_\d> C_{y_k}>C_{y_l}$ if $k > l$ $(j,k,l=1,\cdots,n-1)$.
               \endLR
               From Lemma \ref{L:g2cr}, $C_{y_n}=C_\d$ meets $R$ only along segment $\l_{\d y_n}$. 
               Thus each $C_{y_j}$ runs through $\h_\d^+$ on the side of $\l_{\d y_1}$ 
               $(j=1,\cdots,n-1)$. We have the conclusion from Lemma \ref{L:g2cr}.
               \endR
          \pvcn
               Denote the number of crossings of a flype-tangle $W_x$ by $\dgW_x$, and  
               the number of crossings of $W_x$ which belong to a region $R$ by $\dgW_x|_{R}$. 
               Since $K^*$ is reduced, $W_x$ is trivial if and only if $\dgW_x=0$. 
               Let $p$ be the crossing of $P_x$ which is $n$-adjacent to $\d$. 
               Denote by $U_x$ the region ($\neq P_x$, $Q_x$) which has $x$ but has 
               neither a crossing of $W_x$ nor $\ldq$. Let $u$ (resp. $v$) be the crossing 
               of $U_x$ which is $n$-adjacent (resp. $p$-adjacent) to $x$.

\bgnL{L:g2cx}  Let $x$ be a flype-crossing of $\ldq$ of $K^*$. Then we have the followings.
               \bgnI
               \item[(1)] If $\dgW_x  =   0$, then $C_x=\bar{x}\d$.
               \item[(2)] If $\dgW_x \neq 0$ and $\dgP_x \geq \dgWpx +3$, then $C_x=\bar{x}q\d$. 
               \endI \endLR
               (1) Since $K^*$ is strongly reduced, we have that $p \neq x$.
               Thus applying Lemma \ref{L:g2cc} to $P_x$ and $O_x$, we have the conclusion.
               (2) Since $\dgP_x \geq \dgWpx +3$, we have that $p \neq x$, and 
               thus $C_x$ contains an outside arc $\c_{\bar{x}\d}$ in $P_x$ from Lemma \ref{L:g2cc}. 
               Since $W_x$ is not trivial, $W_x$ has a crossing $x_1$ which is $p$-adjacent to $q$. 
               Let $x_i$ be the crossing of $W_x$ which is $p$-adjacent to $x_{i-1}$ and belongs to 
               $Q_x$ $(i=2,\cdots,n-1)$ so that $x_n=x$. 
               From Lemma \ref{L:g2cc}, $C_{x_1}$ contains an outside arc $\c_{\bar{x}_1 \d}$ in $O_x$.
               If $C_\d >C_{x_2} >C_{x_1}$, then we have that $C_{x_2}=\bar{x}_2q\d\cdots$.
               If $C_{x_2} >C_\d >C_{x_1}$, then $C_\d$ goes into $W_x$ and out from $W_x$ 
               either to $O_x$, to $P_x$ or to $Q_x$. Either case contradicts Lemma \ref{L:g2cr}.
               Then inductively we obtain that $C_{x_n}=C_x=\bar{x}q\d\cdots$. 
               Therefore we can conclude that $C_{x_n}=C_x=\bar{x}q\d$.
               \endR

\bgnL{L:g2cu}  Let $x$ be a flype-crossing of $\ldq$ of $K^*$ and assume that $\dgP_x \geq 
               \dgWpx +4$. If $\dgQ_x=\dgWqx +2$ or $\dgU_x \geq 3$, then $C_u =\bar{u}xq\d$. 
               \endLR
               Since $\dgP_x \geq \dgWpx +4$, we have that $p$ $\neq x,u$. 
               From Lemma \ref{L:g2cx}, we have that $C_x=\bar{x} \d$ if $\dgW_x =0$ and that 
               $C_x=\bar{x}q\d$ if $\dgW_x \neq 0$. Moreover we have that $C_\d > C_u > C_x$ 
               from Lemma \ref{L:g2cc}. Thus we have that $C_u=\bar{u}\d q \cdots$. 
               If $\dgQ_x=\dgWqx +2$, then we have that $v=q$ and thus $C_u$ runs 
               through a side of $\h_x^+$. Hence we have the conclusion.
               Next assume that $\dgQ_x \geq \dgWqx +3$ and $\dgU_x \geq 3$. 
               Then we have that $v \neq u,q$. If $C_\d > C_u > C_v$, then $C_v = u\d q \cdots$, 
               since $C_u=\bar{u}\d q \cdots$. However then $C_v$ admits a non-trivial bridge
               in $Q_x$ or in $U_x$, which contradicts Lemma \ref{L:g2cr}. If $C_\d > C_v > C_u$ 
               or $C_v > C_\d > C_u$, then $C_u$ runs through a side of $\h_x^+$. 
               Hence we have the conclusion.
               \endR

\bgnL{L:g2wx}  Let $x$ be a flype-crossing of $\ldq$ of $K^*$. If $W_x$ is not trivial and 
               does not have a flype-crossing for $\ldq$, then $\h_x$ has a saddle-intersection.
               \endLR
               Since $W_x$ is not trivial, there is a crossing $x_1$ in $W_x$ which is 
               $p$-adjacent to $x$. Then we have that $C_\d > C_x > C_{x_1}$ and that  
               $C_{x_1}$ contains $\c_{\bar{x}_1\d}$ from Lemma \ref{L:g2cc}.
               Moreover since $W_x$ has no flype-crossings for $\ldq$, there 
               exists a crossing $x_2$ in $W_x$ which is $p$-adjacent to $x_1$.
               Then we have that $C_\d > C_x > C_{x_2} > C_{x_1}$ from Lemma \ref{L:g2cr}.
               Therefore $C_{x_2}$ runs through a side of $\h_x^+$. Thus we are done.
               \endR
%
%
%

\bgnP{P:glu2}  If $E$ admits non-equivalent disjoint short cuts, then $K$ has a flyped tongue.
               \endPR
               From Proposition \ref{P:cut2}, $K^*$ is the trivial $2$-component link in 
               prime, strongly reduced almost alternating position. 
               Thus $K^*$ has a flyped tongue from Theorem \ref{T:AAL} (2). 
               Consider $K^*$ with a band-trace $\psi$ of $E$. Note that $\psi$ is
               properly embedded in $R-(E-K)$ for a region $R$. Since the crossings 
               of $K^*$ are preserved by the inverse operation of the extract surgery,
               we use the same notations in $K$ for the crossings of $K^*$.

\bgnM{M:g23a}  Let $x$ be a flype-crossing of $\ldq$ of $K^*$. If $E^* \cap S^+$ has 
               inside arcs $\c_{\d q}$ in $O_x$ and $\c_{xq}$ in $Q_x$, 
               then $K$ admits a flype-component of $\ldq$ with $x$.
               \endMR
               From Lemma \ref{L:g2cr}, 
               the negative curve containing arc $\c_{\d q}$ is $\d q$. 
               Then $\psi$ does not have an end on $\ldq$, and thus $K$ has $\ldq$. 
               Also from Lemma \ref{L:g2cr}, the curve which shares a saddle-intersection
               at $\h_x$ with the curve containing $\c_{xq}$ contains an inside arc
               $\c_{x\d}$ in $P_x$.
               Since $\psi$ does not meet non-boundary arcs, $x$ faces $\d$ (resp. $q$) 
               through $\c_{x \d}$ (resp. $\c_{xq}$) in $S$ with $K$.
               \endR

\bgnM{M:g22f}  If $\ldq$ admits two flype-crossings on one side in $K^*$,
               then $\ldq$ admits a flype-component on the side in $K$.
               \endMR
               Let $x$ and $y$ be flype-crossings of $\ldq$ such that $W_y$ has $x$. 
               We only consider the case when $W_x$ is trivial, 
               since the other case can be shown similarly. 
               We have that $C_x=\bar{x}\d$ and $C_y=\bar{y}q\d$ from Lemma \ref{L:g2cx}. 
               Let $D_x$ (resp. $D_y$) be the disc spanned by $C_x$ (resp. $C_y$) in $S^+-C_\d$.
               Then we can take arcs   
               $\a_x$ in $(P_x \cap D_x) -E$ 
        (resp. $\b_x$ in $(O_x \cap D_x) -E$) 
                                        with ends on $\h_x$ and $\h_\d$ (resp. $\h_q$), and
               $\a_y$ in $(P_x \cap (D_y-D_x)) -E$  
        (resp. $\b_y$ in $(Q_x      -D_y)      -E$)
                                        with ends on $\h_y$ and $\h_\d$ (resp. $\h_q$).
               Note that $\psi$ is properly embedded in $R-(E-K)$ for a region $R$. 
               Thus $y$ faces $\d$ (resp. $q$) through $\a_y$ (resp. $\b_y$) in $K$ if
               $\psi$ is in $D_x -(E-K)$ or in $(D_y \cap Q_x)-(E-K)$, and $x$ faces
               $\d$ (resp. $q$) through $\a_x$ (resp. $\b_x$) in $K$ otherwise.
               \endR

\bgnM{M:g2cd}  Let $x$ be a flype-crossing of $\ldq$ of $K^*$. 
               If $W_x$ is not trivial and $C_\d$ does not run through a side of $\h_x^+$, 
               then $\ldq$ admits a flype-component in $K$ on the same side as $x$.
               \endMR
               If $W_x$ has a flype-crossing for $\ldq$, we are done from Claim \ref{M:g22f}.
               Thus assume otherwise. Then there are positive curves $C$ ($\neq C_\d$) 
               and $C'$ ($\neq C_\d$) which share a saddle-intersection at $\h_x$ such that 
               $C_\d >C >C_x >C'$ from Lemma \ref{L:g2wx} and Lemma \ref{L:g2cc}.
               Therefore $C$ has $\c_{xq}$ in $Q_x$ and $\c_{q\d}$ in $O_x$.
               Hence we are done by Claim \ref{M:g23a}.
               \endR

\bgnM{M:g2d4}  Let $x$ be a flype-crossing of $\ldq$ of $K^*$. If $\dgP_x \geq \dgWpx +4$, 
               then $K$ admits a flype-component of $\ldq$ on the same side with $x$.
               \endMR
               If $\dgU_x=2$, then 
                  $\ldq$ admits a flype-component with $x$ or with $u=v$ from Claim \ref{M:g22f}. 
               If $\dgU_x \geq 3$, then we have that $C_u =\bar{u}xq \d$ from Lemma \ref{L:g2cu}. 
               Thus we are done by Claim \ref{M:g23a}.
               \endR
\pvcn
               Now take a look at a flyped tongue of $K^*$. Let $x_l$ (resp. $x_r$) be a 
               flype-crossing of $\ldq$ of $K^*$ on the left (resp. right) side. We divide the case
               with respect to the degrees of $P_l$, $P_r$ $W_l$ and $W_r$.
               First we have that $\dgP_i \neq 2$, since $K^*$ is strongly reduced ($i=l$ or $r$). 
               Second if $\dgP_l \geq \dgWpl +4$
               and $\dgP_r \geq \dgWpr +4$, then $K$ admits a flyped tongue from Claim \ref{M:g2d4}.
               Third consider the case when $\dgW_l=$ $\dgW_r=0$. 
               Then we may assume that $\dgP_l=3$, since we are done from 
               Claim \ref{M:g2d4} if $\dgP_l \geq 4$ and $\dgP_r \geq 4$. 
               Thus we obtain an alternating diagram $\bar{K}^*$ of the trivial $2$-component link 
               by an untongue move and the Reidemeister move of type II (see Figure \ref{F:ee}). 
               Then $\bar{K}^*$ is disconnected from Theorem \ref{T:AL}. 
               Thus $K^*$ is a diagram of Figure \ref{F:eek} from Theorem \ref{T:AK},
               since $K^*$ is prime. Since $\psi$ is in a region and connects
               different components of $K^*$, $K$ has a flyped tongue.
               Fourth consider the case when $\dgP_i=3$ with $\dgW_i =0$ for $i=l$ or $r$. 
               We are done, since this case is equivalent to the third case.
               Fifth consider the case when $\dgP_l = \dgWpl +3$ and $\dgP_r = \dgWpr +3$. 
               Here we may additionally assume that $\dgW_l \neq 0$ and $\dgW_r \neq 0$ from 
               the above.
               Thus it is sufficient to show that $K$ has a flype-component on the left side
               under the assumption that $C_\d$ runs through a side of $\h_{x_l}^+$
               from Claim \ref{M:g2cd}.
               Then note that $\dgQ_l \geq \dgWql +3$, since otherwise $C_\d$ is not standard. 
               Thus there is a crossing $w$ ($\neq \d$) which is $n$-adjacent to $q$. 
               Then we have that $C_\d>C_w>C_{x_r}$, since $C_{x_r}=\bar{x}_r q \d$ from 
               Lemma \ref{L:g2cx} and thus $C_w=\bar{w}q\d p \cdots$.  
               Therefore the curve sharing a saddle-intersection at $\h_p$ with 
               $C_w$ is $C=px_lq\d$, since $C_{x_l}=\bar{x}_l q \d$ from Lemma \ref{L:g2cx}.
               Hence we are done from Claim \ref{M:g23a}.

\bgnF   {\iclg[scale=.9, bb=157 378 475 493]{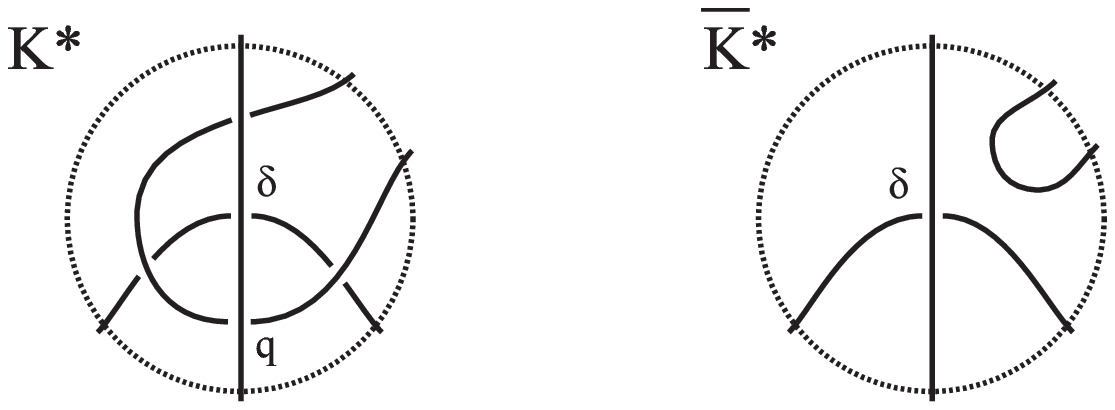}} \captn{} \label{F:ee}\endF 
\bgnF   {\iclg[scale=.9, bb=160 371 487 493]{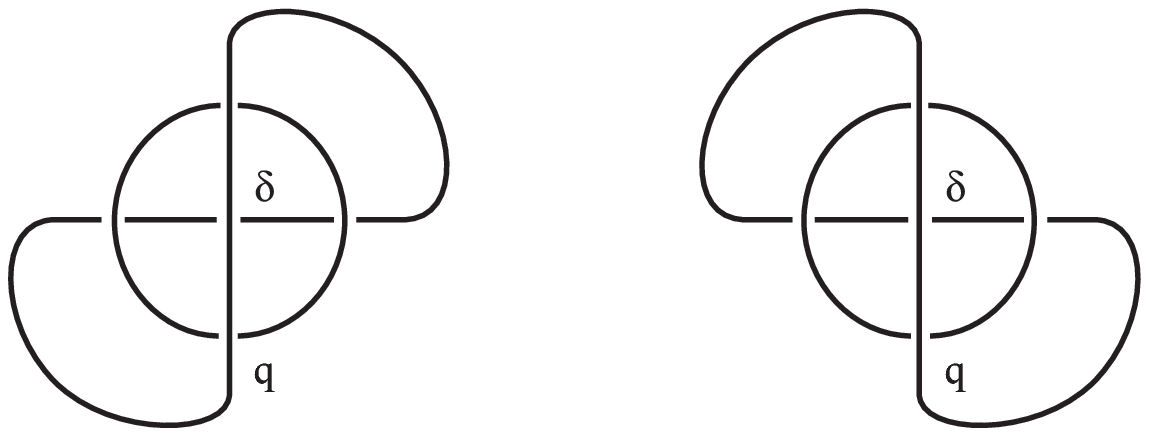}}\captn{} \label{F:eek}\endF

%
%
        \pvcn
               Note that it is impossible to have that $\dgP_l=\dgWpl +2$ and $\dgP_r=\dgWpr +2$, 
               since $\d$ cannot be $n$-adjacent to both $x_l$ and $x_r$. Therefore we are left 
               with the following cases from the symmetry:   
  \pvan
               $\dgP_l=\dgWpl+3$ and $\dgP_r\geq\dgWpr+4$ ($\dgW_l \neq 0$); \\
               $\dgP_l=\dgWpl+2$ and $\dgP_r\geq\dgWpr+4$ ($\dgW_l \neq 0$); and \\
               $\dgP_l=\dgWpl+2$ and $\dgP_r =  \dgWpr+3$ ($\dgW_l \neq 0$ and $\dgW_r \neq 0$).
  \pvcn
               In the first case, let $K^*$ be as the left of Figure \ref{F:g2kj} with 
               regions $A_i$, $B_j$, $I_k$ ($i,j=l,r$ and $k=1, \cdots, 4$) and let $x$ be the 
               crossing of $I_2$ which is $n$-adjacent to $p$. We may assume that $W_l$ has no 
               flype-crossings for $\ldq$, since otherwise we are done by Claim \ref{M:g22f}
               and Claim \ref{M:g2d4}.

\bgnM{M:g2kj}  If $\dgI_1= \dgW_l|_{I_1}+2$, $\dgI_2 =2$ or $\dgI_3 \geq 3$, 
               then $K$ has a flyped tongue.
               \endMR
               From Claim \ref{M:g2d4}, $K$ has a flype-component of $\ldq$ on the right side.
               If $\dgI_1=\dgW_l|_{I_1}+2$, then $C_\d$ does not run through a side 
               of $\h_{x_l}^+$, and thus we are done by Claim \ref{M:g2cd}.
               If $\dg I_2 =2$, then we are done by Claim \ref{M:g22f}. 
               Assume that $\dg I_1\geq \dgW_l|_{I_1}+3$, $\dg I_2\geq 3$ and $\dg I_3 \geq 3$. 
               Then $I_3$ has a crossing $x'$ ($\neq x$) which is $p$-adjacent to $p$. 
               Note that $x \neq x_l$ and $x' \neq x_r$.
               From Claim \ref{M:g2cd} we may assume that $C_\d$ runs through a side of $\h_{x_l}^+$,
               and thus $C_\d = \bar{q}\bar{\d}\bar{p}x_l$. Then we have that $C_\d >C_x >C_{x'}$
               from Lemma \ref{L:g2cr}, 
               since $C_{x'}=\bar{x'}\d q \cdots$ from Lemma \ref{L:g2cc}. Therefore we have that 
               $C_x=\bar{x} p \d q \cdots$, and in fact $C_x$ runs through a side of $\h_p^+$.       
               Hence the positive curve sharing a saddle-intersection with $C_x$ at $\h_p$ 
               is $px_lq\d$, and thus we are done by Claim \ref{M:g23a}.
               \endR
\nwpg  \nidt
               From Claim \ref{M:g2kj} we may assume that $\dg I_1\geq \dgW_l|_{I_1}+3$, 
               $\dg I_2\geq 3$ and $\dg I_3=2$.
               Then we have three subcases: $\dgW_r \neq 0$; $\dgW_r=0$ and $\dg I_4 = 2$; and
               $\dgW_r=0$ and $\dg I_4 \neq 2$.
               In the first subcase, apply flype moves and the untongue move on $K^*$ to have
               another almost alternating diagram $\bK^*$ of the trivial $2$-component link
               with fewer crossings than $K^*$. From given conditions and assumptions, we have
               that $W_l$ is not trivial and has no flype-crossings for $\ldq$, and that 
               $W_r$ is not trivial.
               Let $W_1$ (resp. $W_2$) be the flyped $W_l$ (resp. $W_r$). Since neither $W_l$
               nor $W_r$ is trivial, we have that neither $W_1$ nor $W_2$ is trivial.
               Thus we can see that $\bK^*$ is connected and reduced, and that $\bK^*$ has no
               less than $4$ crossings. Therefore $\bK^*$ is prime from Proposition \ref{P:ctpr}
               and strongly reduced from Proposition \ref{P:g2hl}. Hence $\bK^*$ has a flyped
               tongue from Theorem \ref{T:AAL} (2).
               Let $A_i$, $B_j$ and $H_k$ be mutually distinct regions of $S$ with $\bK^*$ as
               the right of Figure \ref{F:g2kj} ($i,j=1,2$, $k=1,2,3$).
               Let $q'$ be the crossing of $B_2$ which is $n$-adjacent to the dealternator 
               $\d$ of $\bK^*$.
               Since $W_l$ has no flype-crossings for $\ldq$, $W_1$ has no flype-crossings for
               $\l_{\d q'}$, i.e. $B_1$ and $H_1$ do not share a crossing.
               Now each pair of regions $(A_i,B_j)$ $(i,j=1,2)$ can be the core of a 
               flyped tongue of $\bK^*$.
               If $(A_1,B_1)$ is the core of a flyped tongue of $\bK^*$, then $B_2$ shares
               a crossing with the region sharing with $A_1$ a segment which is $n$-adjacent 
               to the crossing of $W_1$ which is $p$-adjacent to the dealternator of $\bK^*$.
               Then the only possibility is $H_1$, contradicting the condition that $W_1$ has
               no flype-crossings. We can be done similarly for $(A_1,B_2)$ and for $(A_2,B_2)$.
               If $(A_2,B_1)$ is the core of a flyped tongue of $\bK^*$, then $B_2$ shares a
               crossing with $H_2$. However then we have that $H_3=B_2$, contradicting that
               $W_2$ is not trivial.
               In the second subcase, the only possibility that $K$ does not have a flype-component 
               of  $\ldq$ with $x_l$ is when band-trace $\psi$ connects segment $\l_{q x_r}$ and  
               the segment $(\neq \l_{x_l p})$ which is $n$-adjacent to $x_l$. 
               However then, $K$ has flype-components of $\l_{\d x_r}$ with $q$ and with $x_l$.
               Now consider the third subcase. If $\dg A_r \neq 4$, i.e. $x$ is not $n$-adjacent 
               to $x_r$, then we obtain a contradiction by applying flype moves and 
               the untongue move on $K^*$ as the first subcase. If $\dg A_r = 4$, then there is
               a crossing $y$ ($\neq p, x_l, x_r$) which is $n$-adjacent to $x$.
               We have that $C_{x_r}=\bar{x}_r q \d$ from Lemma \ref{L:g2cx}, and
               we may assume that $C_\d=\bar{p}\bar{\d}\bar{q}x_l$ from Claim \ref{M:g2cd}.
               Therefore we have that $C_y=\bar{y}x \d q \cdots$ or that $C_y=\bar{y}x p \d q \cdots$.
               In either case there is a positive curve $p\d q x_l$, since $I_3$ has degree $2$.
               Hence we are done from Claim \ref{M:g23a}.

\bgnF   {\iclg[scale=.9, bb=64 406 542 573]{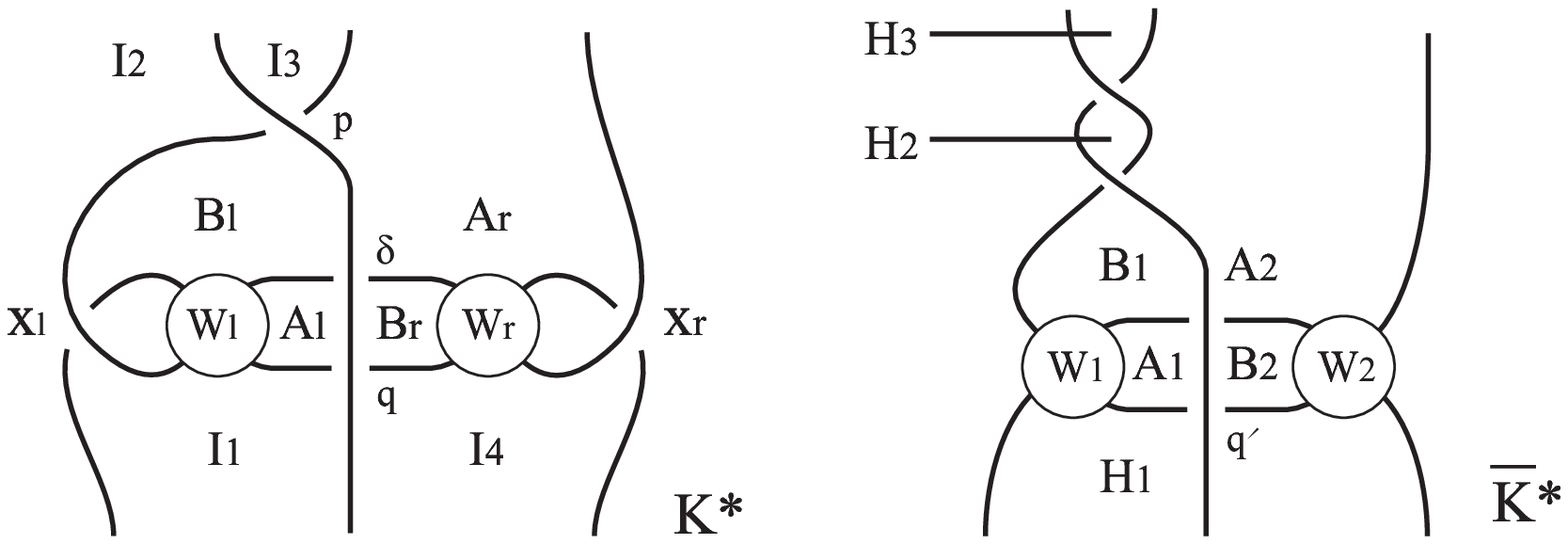}} \captn{} \label{F:g2kj}\endF 

%
%
         \pvcn In the last two cases, we can show that $W_l$ has a flype-crossing for $\ldq$
               by analyzing the diagram after applying flype moves and the untongue move on $K^*$
               as the previous case. However then each case is equivalent to a case which is 
               done before.

\bgnF   {\iclg[scale=.9, bb=78 341 273 496]{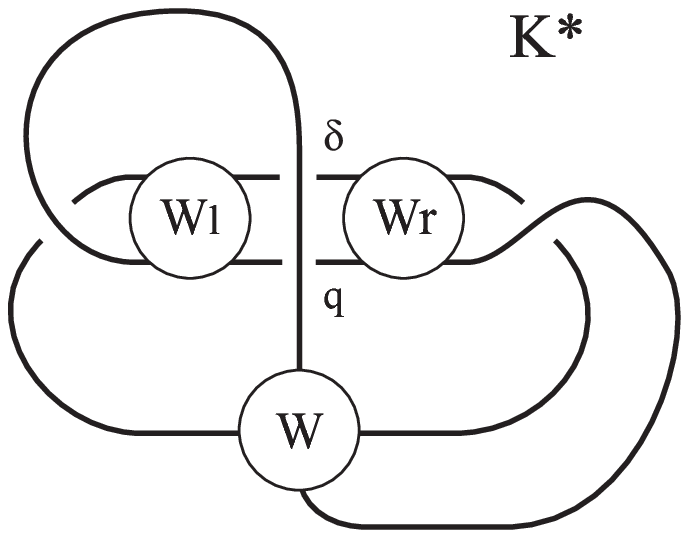}} \captn{} \label{F:g2hl}\endF 

        \endR

\nwpg
\section{Proof of Theorem \ref{T:Main}}                                                
\label{S:Proof}                                                                        

              To prove Theorem \ref{T:AAL} we actually showed the following in \cite{Ts-MPC04}.

\bgnT{T:MPC}  $($\cite{Ts-MPC04}$)$
              Let $L$ be a prime, strongly reduced almost alternating link diagram on $S$.
              Assume that $L$ admits a sphere in $S^3-L$ which is in standard position.
              Then $L$ admits a flyped tongue. \endT 

\noindent{\it Proof of Theorem} \ref{T:Main}.
              Let $K$ be the trivial knot in strongly reduced almost alternating position
              and let $E$ be a spanning disk for $K$ in basic position with minimal complexity.
              Then $K$ is prime from Proposition \ref{P:ctpr} and $E$ is standard position from 
              Proposition \ref{P:MC/ST}. 
              Therefore if $E$ has no short arcs, then $E$ is in special position 
              from Proposition \ref{P:nosa}.
              Then we can take a neighborhood of $E$ whose boundary is in standard position
              from Proposition \ref{P:SP/ST}.
              Thus $K$ has a flyped tongue from Theorem \ref{T:MPC}. Next assume that $E$ has
              a short arc. If $E$ has a short arc of type I, then $K$ has a flyped tongue from
              Proposition \ref{P:glu1}. If $E$ has a short arc of type II or III, then $K$ has 
              a flyped tongue from Proposition \ref{P:a2sc} and Proposition \ref{P:glu2}.
              \End




\begin{thebibliography}{99}

\bibitem[Ad1]{Ad-TPAP92}       C. Adams,  J. Brock,  J. Bugbee,        T. Comar,  
                               K. Faigin, A. Huston, A. Joseph  and    D. Pesikoff, 
                          {\em Almost alternating links}, 
                               Topology Appl. 
                          {\bf 46} (1992), 151--165.

\bibitem[Ad2]{Ad-93}           C. Adams,  C. Arthur, D. Bruneau, T. Graber, 
                               J. Kucera, P. Vongsinsirikul and  T. Welsh, 
                          {\em The reduction of almost alternating links and knots}, 
                               preprint (1993).

\bibitem[Ad3]{Ad-94}           C. Adams, 
                          {\em The knot book; An elementary introduction to 
                               the mathematical theory of knots}, 
                               W.H. Freeman and Company, New York (1994). 

\bibitem[Ba]{Ba-MA30}          C. Bankwitz, 
                          {\em $\ddot{U}$ber die Torsionszahlen der alternierenden Knoten}, 
                               Math. Ann. 
                          {\bf 103} (1930), 145--161.

\bibitem[BM]{BM-TAMS92}        J.S. Birman and W.W. Menasco,
                          {\em Studying links via closed braids ${\rm V}$: The unlink}, 
                               Trans. AMS 
                          {\bf 329} (1992), 585--606.

\bibitem[BZ]{BZ-85}            G. Burde and H. Zieschang,
                          {\em Knots}, 
                               de Gruyter, 1985.

\bibitem[Crm]{Crm-JLMC59}      P.R. Cromwell,
                          {\em Homogeneous links}, 
                               J. London Math. Soc. (2)
                          {\bf 39} (1989), 535--552.

\bibitem[Crw]{Crw-ANN59}       R.H. Crowell,
                          {\em Genus of alternating link types}, 
                               Ann. of Math. (2)
                          {\bf 69} (1959), 258--275.

\bibitem[GL]{GL-06}            C. McA. Gordon and J. Luecke,
                          {\em Knots with unknotting number $1$ and essential Conway spheres}, 
                               preprint, arXiv:math.GT/0601265. 

\bibitem[Hi]{Hi-TAP00}         M. Hirasawa, 
                          {\em Triviality and splittability of special almost-alternating 
                               diagrams via canonical Seifert surfaces}, 
                               Topology Appl. 
                          {\bf 102} (2000), 89--100.

\bibitem[Ko]{Ko-PAMS91}         P. Kohn,
                          {\em Two-bridge links with unlinking number one}, 
                               Proc. Amer. Math. Soc.
                          {\bf 113} (1991), 1135--1147.


\bibitem[Me]{WM-TPL84}         W.W. Menasco, 
                          {\em Closed incompressible surfaces in alternating knot 
                               and link complements}, 
                               Topology 
                          {\bf 23} (1984), 37--44.

\bibitem[MT1]{MT-MPC91}        W.W. Menasco and M.B. Thistlethwaite, 
                          {\em A geometric proof that alternating knots are non-trivial}, 
                               Math. Proc. Camb. Philos. Soc. 
                          {\bf 109} (1991), 425--431.

\bibitem[MT2]{MT-JRAM92}        W.W. Menasco and M.B. Thistlethwaite, 
                          {\em Surfaces with boundary in alternating knot exteriors}, 
                               J. Reine Angew. Math. 
                          {\bf 426} (1992), 47--65.

\bibitem[Mu]{Mu-JMSJ58}        K. Murasugi, 
                          {\em On the genus of the alternating knot {\rm I}, {\rm II}}, 
                               J. Math. Soc. Japan
                          {\bf 10} (1958), 94--105, 235--248.

\bibitem[St]{St-CM04}          A. Stoimenow,
                          {\em Gau${\rm \beta}$ diagram sums on almost positive knots},
                               Compositio Math. 
                          {\bf 137} (2004) 228--254.  

\bibitem[Ts]{Ts-MPC04}         T. Tsukamoto,
                          {\em A criterion for almost alternating links to be non-splittable},
                               Math. Proc. Camb. Philos. Soc. 
                          {\bf 137} (2004) 109--133.  


                                        \end{thebibliography}
                                            \end{document}